\documentclass[a4paper,10pt]{article}

\usepackage{graphicx}
\usepackage{amsmath}
\usepackage[authoryear]{natbib}
\usepackage{amsthm}
\usepackage{microtype}
\usepackage{bm}
\usepackage[hmargin=2.5cm,vmargin={3.5cm,4cm}]{geometry}
\usepackage[bitstream-charter]{mathdesign}
\usepackage[final]{showkeys}

\numberwithin{equation}{section}
\allowdisplaybreaks

\newtheorem{Theorem}{Theorem}[section]
\newtheorem{Remark}{Remark}[section]

\newtheorem{Lemma}{Lemma}[section]

\def\sgn{\mathop{\rm sgn}}

\author{$\text{Enzo Orsingher}_1$, $\text{Federico Polito}_2$\\
	\footnotesize (1) -- Dipartimento di Scienze Statistiche,
	``Sapienza'' Universit\`a di Roma\\
	\footnotesize Piazzale Aldo Moro 5, 00185 Rome, Italy.\\
	\footnotesize Tel: +39-06-49910585, fax: +39-06-4959241\\
	\footnotesize Email address: enzo.orsingher@uniroma1.it (corresponding author)\\
	\footnotesize (2) -- Dipartimento di Matematica, Universit\`a degli studi di Roma
	``Tor Vergata''\\
	\footnotesize Via della Ricerca Scientifica 1, 00133, Rome, Italy.\\
	\footnotesize Tel: +39-06-2020568, fax: +39-06-20434631\\
	\footnotesize Email address: polito@nestor.uniroma2.it
	}

\title{Compositions, Random Sums and Continued Random Fractions of Poisson and Fractional Poisson Processes}

\begin{document}

	\maketitle

	\begin{abstract}
		
		\noindent In this paper we consider the relation between random sums and compositions of different processes.
		In particular, for independent Poisson processes $N_\alpha(t)$, $N_\beta(t)$, $t>0$, we show that
		$N_\alpha(N_\beta(t)) \overset{\text{d}}{=} \sum_{j=1}^{N_\beta(t)} X_j$, where the $X_j$s are
		Poisson random variables. We present a series of similar cases, the most general of which is the one in which
		the outer process is Poisson and the inner one is a nonlinear fractional birth process. We highlight
		generalisations of these results where the external process is infinitely divisible.
		A section of the paper concerns compositions of the form $N_\alpha(\tau_k^\nu)$, $\nu \in (0,1]$,
		where $\tau_k^\nu$ is the inverse of the fractional Poisson process, and we show how these compositions can be
		represented as random sums.
		Furthermore we study compositions of the form $\Theta(N(t))$, $t>0$, which can be represented as
		random products.
		The last section is devoted to studying continued fractions of Cauchy random variables
		with a Poisson number of levels. We evaluate the exact distribution and derive the scale parameter
		in terms of ratios of Fibonacci numbers.
		
		\medskip
		
		
		\medskip

	\end{abstract}

	\section{Introduction}
	
		Publications in the field of probability have devoted considerable attention to
		compositions of
		different processes, as e.g.\ Brownian motions, fractional Brownian motions and telegraph processes.
		Also, more general cases as stable processes, combined in different ways have been investigated
		and the p.d.e.\ connections analysed.
		In the present paper we focus our attention on the composition of point processes, e.g.\
		Poisson processes, fractional Poisson processes, and others.
		For independent homogeneous Poisson processes $N_\alpha(t)$, $N_\beta(t)$, $t>0$, we are
		able to show that $N_\alpha(N_\beta(t))$ has a remarkable connection, with respect to
		distributional properties, with random sums,
		that is, we prove that
		\begin{align}
			\label{pa}
			\hat{N}(t) = N_\alpha(N_\beta(t)) \overset{\text{d}}{=} \sum_{j=1}^{N_\beta(t)} X_j.
		\end{align}
		The first part of the paper is devoted to the presentation of similar results involving more general processes
		as e.g.\ nonlinear fractional birth processes $\mathcal{Y}^\nu(t)$, $t>0$, $\nu \in (0,1]$, so as to obtain
		\begin{align}
			\label{pa2}
			N_\alpha(\mathcal{Y}^\nu(t)) \overset{\text{d}}{=} \sum_{j=1}^{\mathcal{Y}^\nu(t)} X_j,
		\end{align}
		where the $X_j$s appearing in \eqref{pa} and \eqref{pa2} are independent Poisson random variables
		of parameter $\lambda_\alpha$.
		A more general result is obtained when the external process is replaced by a process $\Theta$ possessing
		infinitely divisible distribution. In this case, we are able to show that
		\begin{align}
			\Theta(N_\beta(t)) \overset{\text{d}}{=} \sum_{j=1}^{N_\beta(t)} \xi_j,
		\end{align}
		where the random variables $\xi_j$s are the components of the infinitely divisible random variable
		$\Theta(1)$, and $N_\beta(t)$, $t>0$, is a homogeneous Poisson process.
		The representation \eqref{pa} is remarkable in that it results in the explicit
		law of a random sum which is rarely possible in general. The distribution of \eqref{pa}
		can be given as
		\begin{align}
			\text{Pr} \{ N_\alpha(N_\beta(t)) = k \} & =
			\frac{\lambda_\alpha^k}{k!} e^{-\lambda_\beta t} \sum_{r=0}^\infty \frac{e^{-\lambda_\alpha r}
			r^k(\lambda_\beta t)^r}{r!} \\
			& = \frac{\lambda_\alpha^k}{k!} e^{-\lambda_\beta t (1-e^{-\lambda_\alpha})}
			\mathfrak{B}_k \left( \lambda_\beta t e^{-\lambda_\alpha} \right), \notag
		\end{align}
		where
		\begin{align}
			\mathfrak{B}_k(x) = e^{-x} \sum_{r=0}^\infty \frac{r^k x^r}{r!},
		\end{align}
		are the so-called Bell polynomials, and $\lambda_\alpha$, $\lambda_\beta$, are the parameters of,
		respectively, $N_\alpha(t)$ and $N_\beta(t)$.
		The composition \eqref{pa} produces a process with linearly increasing mean value and variance:
		\begin{align}
			\begin{cases}
				\mathbb{E}N_\alpha(N_\beta(t)) = \lambda_\alpha \lambda_\beta t, \\
				\mathbb{V}\text{ar} N_\alpha(N_\beta(t)) = \lambda_\alpha (\lambda_\alpha + 1)\lambda_\beta t.
			\end{cases}
		\end{align}
		
		For the iterated Poisson process we find that the first-passage time
		\begin{align}
			T_k = \inf (s\colon N_\alpha(N_\beta(s))=k ), \qquad k \geq 1,
		\end{align}
		has distribution
		\begin{align}
			\text{Pr} \{ T_k \in \mathrm ds \} = \mathrm ds \, \lambda_\beta e^{-\lambda_\alpha} e^{-\lambda_\beta s}
			\frac{\lambda_\alpha^k}{k!} \sum_{j=0}^\infty e^{-\lambda_\alpha j} \left[ (j+1)^k-j^k \right]
			\frac{(\lambda_\beta s)^j}{j!}, \qquad s<0.
		\end{align}
		Furthermore, $\text{Pr} \{ T_k < \infty \} < 1$, for all $k\geq 1$, so that there is a positive probability
		of never hitting the level $k$ because the iterated Poisson process can take jumps of arbitrary
		integer-valued size.	
		
		Finally, we note that the iterated Poisson process produces a Galton--Watson process (continuous in time).

		In the case in which $N_\beta(t)$ is replaced by a non-homogeneous Poisson process $\mathfrak{N}(t)$,
		with rate $\lambda(t)$, $t>0$, we still have a representation of the composition $N_\alpha(\mathfrak{N}(t))$
		as a random sum, i.e.\
		\begin{align}
			\label{teaway}
			N_\alpha(\mathfrak{N}(t)) \overset{\text{d}}{=} \sum_{j=1}^{\mathfrak{N}(t)} X_j,
		\end{align}
		where the $X_j$s are independent Poisson random variables of parameter $\lambda_\alpha$.
		Since the probability generating function of \eqref{teaway} reads
		\begin{align}
			\mathbb{E} u^{N_\alpha(\mathfrak{N}(t))} = e^{\int_0^t\lambda(w)\mathrm dw\left[ e^{\lambda_\alpha
			(u-1)} -1 \right]}
		\end{align}
		we have also that
		\begin{align}
			\label{limton}
			\begin{cases}
				\mathbb{E} N_\alpha(\mathfrak(t)) = \lambda_\alpha \int_0^t \lambda(w) \mathrm dw, \\
				\mathbb{V}\text{ar} = \lambda_\alpha(\lambda_\alpha+1) \int_0^t \lambda(w) \mathrm dw.
			\end{cases}
		\end{align}
		The process \eqref{pa} has non-decreasing sample paths with jumps of integer-valued
		size and thus differs and extends the classical homogeneous Poisson process.
		We note that, for $\lambda_\alpha \neq \lambda_\beta$, we have that
		\begin{align}
			\text{Pr} \left\{ N_\alpha(N_\beta(t)) \neq N_\beta(N_\alpha(t)) \right\} = 1,
		\end{align}
		as can be inferred from the structure of the probability generating function.
		The interchange of the non-homogeneous Poisson process $\mathfrak{N}(t)$, with the homogeneous
		one produces a composed process
		\begin{align}
			\label{milton}
			\mathfrak{N}(N_\alpha(t)),
		\end{align}
		whose probability generating function can be written as
		\begin{align}
			\mathbb{E}u^{\mathfrak{N}(N_\alpha(t))} = \sum_{r=0}^\infty e^{(u-1)\int_0^r\lambda(w) \mathrm dw}
			\frac{(\lambda_\alpha t)^r}{r!} e^{-\lambda_\alpha t}.
		\end{align}
		Even the mean value of \eqref{milton} differs from \eqref{limton} since
		\begin{align}
			\mathbb{E}\mathfrak{N}(N_\alpha(t)) & = e^{-\lambda_\alpha t} \sum_{r=1}^\infty \left[ 
			\int_0^r \lambda(w) \mathrm dw \right] \frac{(\lambda_\alpha t)^r}{r!}
			= \sum_{j=1}^\infty \left[ \int_{j-1}^j \lambda(w) \mathrm dw \right]
			\text{Pr} \{ N_\alpha(t) \geq j \}.
		\end{align}
		
		The third section deals with the inverse of the fractional Poisson process $N^\nu(t)$, $t>0$.
		The process $N^\nu(t)$ can be viewed as a renewal process with Mittag--Leffler distributed intertimes
		and the following probability distribution (see \citet{mainardi,beg})
		\begin{align}
			\text{Pr} \{ N^\nu(t) = m \} & = \sum_{j=m}^\infty (-1)^{j-m} \binom{j}{m} \frac{(\lambda_\beta
			t^\nu)^j}{\Gamma(\nu j +1)} \\
			& = (\lambda_\beta t^\nu)^m
			\sum_{j=0}^\infty \binom{j+m}{j} \frac{(-\lambda_\beta t^\nu)^j}{\Gamma \left(\nu(m+j)+1\right)},
			\qquad t>0, \: \nu \in (0,1]. \notag
		\end{align}
		The inverse of the fractional Poisson process is defined as
		\begin{align}
			\tau^\nu_k = \inf (t \colon N^\nu(t) = k), \qquad k \geq 1, \: \nu \in (0,1],
		\end{align}
		with distribution
		\begin{align}
			\text{Pr} \{ \tau_k^\nu \in \mathrm ds \} / \mathrm ds =
			\lambda_\beta^k \sum_{j=0}^\infty \binom{-k}{j} \lambda_\beta^j \frac{s^{\nu(k+j)-1}}{\Gamma(\nu(k+j))}
		\end{align}
		and moment generating function
		\begin{align}
			\mathbb{E} e^{-\mu \tau_k^\nu} = \left( \frac{\mu^\nu}{\lambda_\beta} +1 \right)^{-k}.
		\end{align}
		It should be noted that for $\nu = 1$ this coincides with the Erlang process.
		The composed process $N_\alpha(\tau_k^\nu)$ has probability generating function
		\begin{align}
			\label{stir}
			\mathbb{E}u^{N_\alpha(\tau_k^\nu)} = \left[ 1+(1-u)^\nu \frac{\lambda_\alpha^\nu}{\lambda_\beta} \right]^{-k}
		\end{align}
		which suggests the following interesting representation:
		\begin{align}
			N_\alpha(\tau_k^\nu) \overset{\text{d}}{=} \sum_{j=1}^k \xi_j,
		\end{align}
		where the independent random variables $\xi_j$, $1\leq j \leq k$,
		are discrete Mittag--Leffler (see \citet{pillai}).
		From \eqref{stir}, we can infer that $N_\alpha(\tau_k^\nu)$ has Linnik distribution and, for $\nu=1$, this
		coincides with the negative binomial distribution having parameters $k$ and $\lambda_\alpha/(\lambda_\alpha+
		\lambda_\beta)$.
		For the special case $\nu=1$, we also have the following representation of $N_\alpha(\tau_k^\nu)$:
		\begin{align}
			N_\alpha(\tau_k^\nu) \overset{\text{d}}{=} X_1 + \dots + X_N,
		\end{align}
		where $N$ is a Poisson random variable of parameter
		$\mu = \log ((\lambda_\alpha+\lambda_\beta)/\lambda_\beta)^k$,
		and the $X_j$s are independent and possess logarithmic distribution of parameter
		$\eta = \lambda_\alpha/(\lambda_\alpha + \lambda_\beta)$.
		Furthermore, the case $N_\alpha(\phi_k^\nu)$ where the inner process is the inverse $\phi_k^\nu$ of a
		fractional linear pure birth process $Y_\beta^\nu(t)$, is examined.
		
		In particular, we show that
		\begin{align}
			\mathbb{E}u^{N_\alpha(Y_\beta^\nu(t))} =
			k! \frac{\Gamma\left( \frac{\lambda_\alpha^\nu(1-u)^\nu}{\lambda_\beta} +1 \right)}{
			\Gamma\left( \frac{\lambda_\alpha^\nu(1-u)^\nu}{\lambda_\beta} +1+k \right)}, \qquad |u|<1.
		\end{align}
		
		The final section of this paper deals with random products of
		non negative i.i.d.\ random variables of the form
		\begin{align}
			N_\pi (t) = \prod_{i=1}^{N(t)} X_j, \qquad t>0,
		\end{align}
		where $N(t)$ is a homogeneous Poisson process. We show that the Mellin transform of
		$N_\pi(t)$ is
		\begin{align}
			\mathbb{E} \left( N_\pi(t) \right)^{\eta -1} = e^{\lambda t \left[ \mathbb{E}
			X^{\eta -1} -1 \right]}.
		\end{align}
		Furthermore, we can evaluate the covariance function of the process $N_\pi(t)$ which reads
		\begin{align}
			\mathbb{C}\text{ov} \left[ N_\pi(t), N_\pi(s) \right]
			= e^{\lambda t (\mathbb{E} X -1)} \left( e^{\lambda s \mathbb{E}X(X-1)} -
			e^{\lambda s(\mathbb{E}X -1)} \right).
		\end{align}
		Clearly, the random products can be viewed as compositions of the form $\Xi(N(t))$, where
		\begin{align}
			\Xi(k)=
			\begin{cases}
				1, & k=0, \\
				\prod_{j=1}^k X_j = e^{\sum_{j=1}^k \log X_j}, & k > 1.
			\end{cases}
		\end{align}
		Finally, we consider continuous fractions with a random number of levels:
		\begin{align}
			\label{contf}
			[ X_1;X_2,\dots,X_{N(t)} ] = X_1 + \frac{1}{X_2 + \frac{1}{\ddots + X_{N(t)-1} + \frac{1}{X_{N(t)}}}},
		\end{align}
		where the $X_j$s are independent standard Cauchy random variables and $N(t)$ is a homogeneous
		Poisson process. We show that the conditional distribution
		of \eqref{contf} is Cauchy in which the scale parameter equals $F_{n+1}/F_n$, where $F_n$ are Fibonacci
		numbers.
		This permits us to give a stochastic representation of \eqref{contf} in the form
		\begin{align}
			[X_1;X_2, \dots, X_{N(t)}] \overset{\text{d}}{=} \sum_{j=1}^{F_{N(t) +1}} Y_{j,N(t)},
		\end{align}
		where $Y_{j,N(t)}$ are Cauchy random variables with scale parameter equal to $1/F_{N(t)}$.
			
	\section{Composition of Poisson processes with different point processes}
	
		In this section we examine the following compositions:
		
		\begin{enumerate}
			\item $N_\alpha(N_\beta(t))$, where $N_\alpha$ and $N_\beta$ are independent Poisson processes;
			\item $N_\alpha(\mathfrak{N}(t))$, where $\mathfrak{N}(t)$ is a non-homogeneous Poisson process
				with rate $\lambda(t)$, $t>0$;
			\item $N_\alpha(\mathcal{Y}^\nu(t))$, where $\mathcal{Y}^\nu(t)$ is a nonlinear fractional birth
				process with birth rates $\lambda_1,\dots,\lambda_k,\dots$, and $0<\nu\leq 1$;
			\item $N_\alpha(Y^\nu(t))$, where $Y^\nu(t)$ is a linear fractional birth process. This is a special
				case of the preceding point.
			\item $N_\alpha(\mathcal{F}(\mathcal{B}))$, where $\mathcal{F}(\mathcal{B})$ is a
				Poisson field.
			\item $\Theta(N_\beta(t))$, where $\Theta(1)$ is an infinitely divisible random variable
				and $N(t)$ is an arbitrary point process. 
		\end{enumerate}
		We will establish some distributional relations between these composed processes and random sums.
		
		\subsection{Iterated Poisson process}
		
\begin{figure}
				\centering
				\includegraphics[scale=0.5]{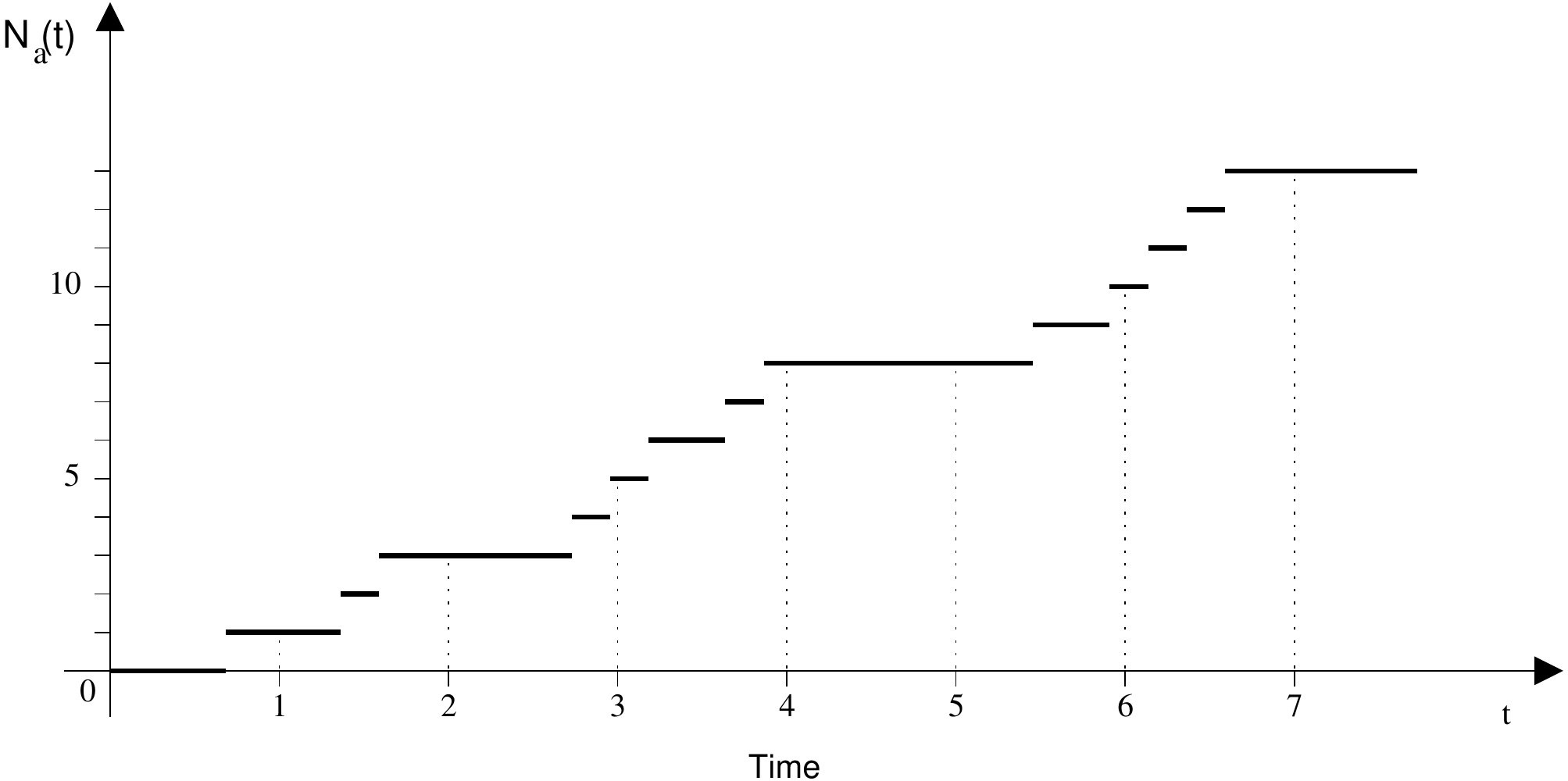}
				\caption{\label{figpath1}A realisation of the external Poisson process
				$N_\alpha(t)$, $t>0$.}
			\end{figure}

\begin{figure}
				\centering
				\includegraphics[scale=0.5]{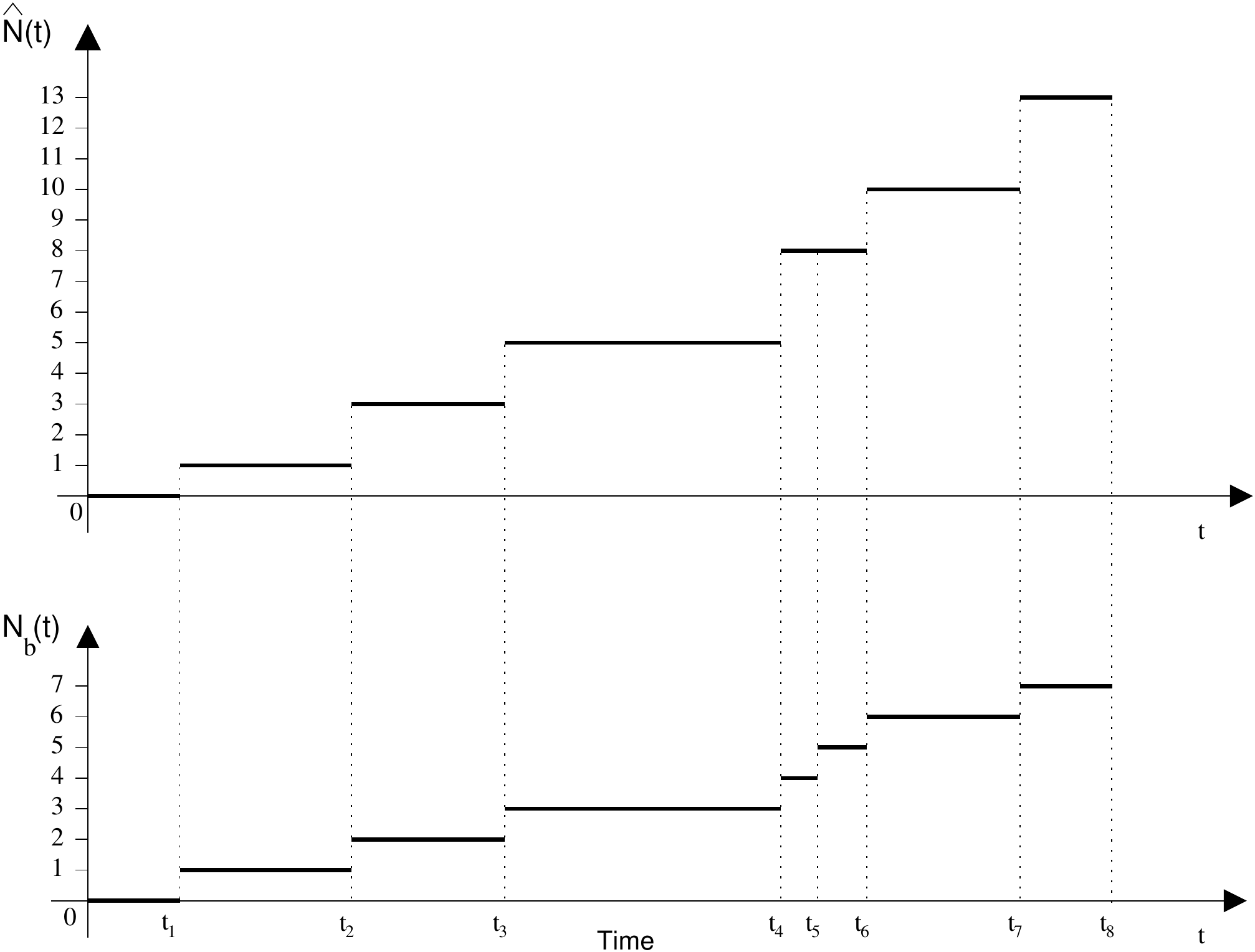}
				\caption{\label{figpath2}A path of the process
				$\hat{N}(t)=N_\alpha(N_\beta(t))$, $t>0$, together with its relation for a specific
					path of the internal process $N_\beta(t)$, $t>0$.}
			\end{figure}

			We start our analysis by considering the iterated Poisson process $\hat{N}(t)=N_\alpha(N_\beta(t))$,
			$t>0$. The sample paths of $\hat{N}(t)$ are non-decreasing with jumps of arbitrary integer-valued size.
			In the figures \ref{figpath1} and \ref{figpath2}, we give the trajectories of $N_\alpha(t)$ and $N_\beta(t)$
			separately and then the sample path of $\hat{N}(t)$ obtained by their composition.
	
			We note that the iterated Poisson process $\hat{N}(t)$ jumps at the occurence of events of the inner process
			$N_\beta(t)$. Thus, if the rate of $N_\beta(t)$ is large, then $\hat{N}(t)$ has rapidly increasing
			trajectories for large $\lambda_\alpha$. However high values of $\lambda_\alpha$ and low
			levels of $\lambda_\beta$ can produce contradictory results and thus compensate each other.
			
			\begin{Theorem}
				The distribution of $\hat{N}(t)=N_\alpha(N_\beta(t))$ reads
				\begin{align}
					\label{elena}
					\text{Pr} \{ \hat{N}(t)=k \} & =
					\frac{\lambda_\alpha^k}{k!} e^{-\lambda_\beta t} \sum_{r=0}^\infty \frac{e^{-\lambda_\alpha r}
					r^k(\lambda_\beta t)^r}{r!} \\
					& = \frac{\lambda_\alpha^k}{k!} e^{-\lambda_\beta t (1-e^{-\lambda_\alpha})}
					\mathfrak{B}_k \left( \lambda_\beta t e^{-\lambda_\alpha} \right), \qquad k \geq 0, \: t>0, \notag
				\end{align}
				where
				\begin{align}
					\mathfrak{B}_k(x) = e^{-x} \sum_{r=0}^\infty \frac{r^k x^r}{r!},
				\end{align}
				is the $k$th order Bell polynomial \citep{khristo}). The probability generating function of
				$\hat{N}(t)$ has the form
				\begin{align}
					\label{mah}
					\mathbb{E}u^{\hat{N}(t)} = e^{\lambda_\beta t \left( e^{\lambda_\alpha(u-1)}-1 \right)},
					\qquad |u|\leq 1.
				\end{align}
				\begin{proof}
					\begin{align}
						\text{Pr} \{ \hat{N}(t)=k \} & = \sum_{r=0}^\infty \text{Pr} \{ N_\alpha(r)=k \}
						\text{Pr} \{ N_\beta(t)=r \} \\
						& = \sum_{r=0}^\infty e^{-\lambda_\alpha r} \frac{(\lambda_\alpha r)^k}{k!}
						e^{-\lambda_\beta t} \frac{(\lambda_\beta t)^r}{r!}. \notag
					\end{align}
					Furthermore
					\begin{align}
						\mathbb{E}u^{\hat{N}(t)} & = \sum_{k=0}^\infty u^k \sum_{r=0}^\infty
						\frac{e^{-\lambda_\alpha r}(\lambda_\alpha r)^k}{k!}
						\cdot \frac{e^{-\lambda_\beta t}(\lambda_\beta t)^r}{r!} \\
						& = \sum_{r=0}^\infty e^{-\lambda_\alpha r} e^{u \lambda_\alpha r}
						e^{\lambda_\beta t} \frac{(\lambda_\beta t)^r}{r!}, \notag
					\end{align}
					and thus \eqref{mah} emerges.
				\end{proof}
			\end{Theorem}

\begin{figure}
				\centering
				\includegraphics[scale=1]{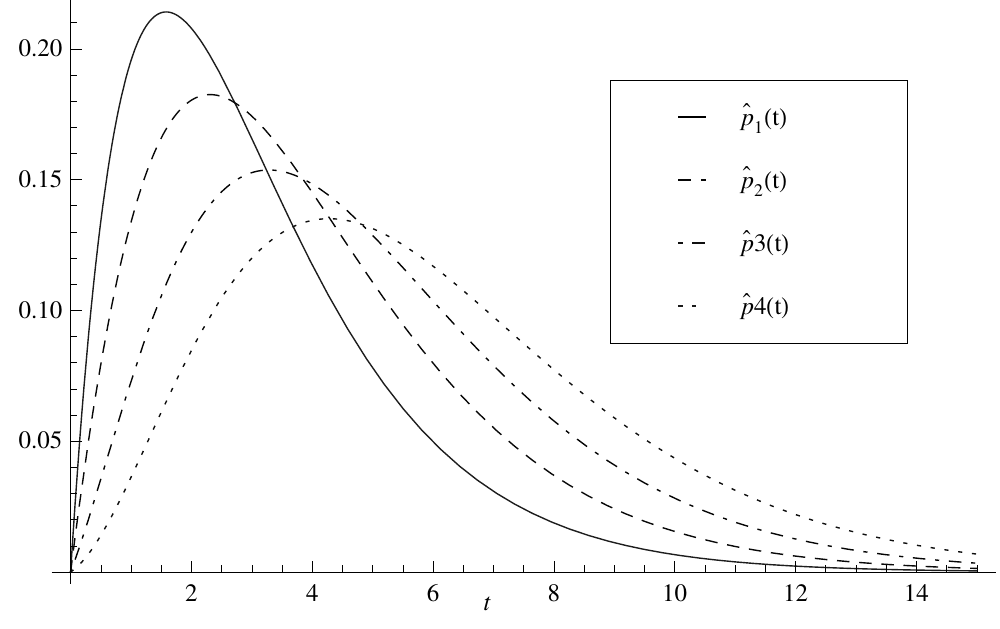}
				\caption{\label{figprob}The first four state probabilities (Iterated Poisson).
					The parameters are $\lambda_\alpha=1$, $\lambda_\beta=1$.}
			\end{figure}

			In Figure \ref{figprob} we present the first four state probabilities as a function of time $t$
			for the iterated Poisson
			process in which $\lambda_\alpha=\lambda_\beta=1$.
		
			\begin{Theorem}
				The following equality in distribution holds:
				\begin{align}
					\label{mar}
					N_\alpha(N_\beta(t)) \overset{\text{d}}{=} \sum_{j=1}^{N_\beta(t)} X_j,
				\end{align}
				where the $X_j$s are i.i.d.\ Poisson random variables of parameter $\lambda_\alpha$.
				\begin{proof}
					The probability generating function of the random sum \eqref{mar} is
					\begin{align}
						\mathbb{E}u^{\sum_{j=1}^{N_\beta(t)} X_j} & = \sum_{k=0}^\infty
						\left( \mathbb{E}u^X \right)^k \text{Pr} \{ N_\beta(t)=k \}
						= e^{-\lambda_\beta t + \lambda_\beta t \mathbb{E}u^X}
						= e^{\lambda_\beta t \left(e^{\lambda_\alpha(u-1)}-1  \right)},
					\end{align}
					and this coincides with \eqref{mah}.
				\end{proof}
			\end{Theorem}
		
			\begin{Remark}
				For the iterated Poisson process we can write
				\begin{align}
					\label{beee}
					\mathbb{E}\hat{N} (t) = \mathbb{E} N_\alpha(N_\beta(t)) = \lambda_\alpha \lambda_\beta t,
				\end{align}
				\begin{align}
					\label{veee}
					\mathbb{V}\text{ar} \hat{N}(t) = \lambda_\alpha (1+\lambda_\alpha) \lambda_\beta t.
				\end{align}
				We can obtain \eqref{beee} directly, by means of the probability generating function,
				and also by applying Wald's formula for random sums.
				\begin{align}
					\mathbb{E} N_\alpha(N_\beta(t)) & = \sum_{k=0}^\infty k \sum_{r=0}^\infty
					e^{-\lambda_\alpha r} \frac{(\lambda_\alpha r)^k}{k!} e^{-\lambda_\beta t} \frac{(\lambda_\beta t)^r}{
					r!} \\
					& = e^{-\lambda_\beta t} \sum_{r=0}^\infty e^{-\lambda_\beta r} \frac{(\lambda_\beta t)^r}{
					r!} \sum_{k=1}^\infty \frac{(\lambda_\alpha r)^k}{(k-1)!} \notag \\
					& = e^{-\lambda_\beta t} \sum_{r=1}^\infty \lambda_\alpha^r \frac{(\lambda_\beta t)^r}{r!}
					= \lambda_\alpha \lambda_\beta t. \notag
				\end{align}
				Alternatively
				\begin{align}
					\mathbb{E}N_\alpha(N_\beta(t)) = \left. \frac{\mathrm d}{\mathrm du} \mathbb{E}u^{\hat{N}(t)}
					\right|_{u=1}
					= \left. \lambda_\alpha \lambda_\beta t \left( e^{\lambda_\alpha(u-1)} \right)
					e^{\lambda_\beta t\left(e^{\lambda_\alpha(u-1)-1}\right)} \right|_{u=1} =
					\lambda_\alpha \lambda_\beta t
				\end{align}
				and, by applying the Wald's formula
				\begin{align}
					\mathbb{E}X \mathbb{E}N_\beta (t) = \mathbb{E} \hat{N}(t) = \lambda_\alpha \lambda_\beta t.
				\end{align}
				For the second moment, analogously we obtain
				\begin{align}
					\mathbb{E}\hat{N}^2 (t) & = \sum_{k=0}^\infty k^2 \sum_{r=0}^\infty e^{-\lambda_\alpha r}
					\frac{(\lambda_\alpha r)^k}{k!} e^{-\lambda_\beta t} \frac{(\lambda_\beta t)^r}{r!} \\
					& = e^{-\lambda_\beta t} \sum_{r=0}^\infty e^{-\lambda_\alpha r} \frac{(\lambda_\beta t)^r}{r!}
					\sum_{k=0}^\infty \frac{k^2}{k!} (\lambda_\alpha r)^k \notag \\
					& = e^{-\lambda_\beta t} \sum_{r=0}^\infty e^{-\lambda_\alpha r} \frac{(\lambda_\beta t)^r}{r!}
					\sum_{k=0}^\infty \frac{k+1}{k!} (\lambda_\alpha r)^{k+1} \notag \\
					& = e^{-\lambda_\beta t} \sum_{r=0}^\infty e^{-\lambda_\alpha r} \frac{(\lambda_\beta t)^r}{r!}
					\left[ \sum_{k=0}^\infty \frac{(\lambda_\alpha r)^{k+2}}{k!} + \sum_{k=0}^\infty
					\frac{(\lambda_\alpha r)^{k+1}}{k!} \right] \notag \\
					& = e^{-\lambda_\beta t} \sum_{r=0}^\infty \frac{(\lambda_\beta t)^r}{r!}
					\left[ (\lambda_\alpha r)^2 + (\lambda_\alpha r) \right] \notag \\
					& = e^{-\lambda_\beta t} \left[ \sum_{r=0}^\infty \frac{(\lambda_\beta t)^{r+1}}{r!}
					\lambda_\alpha^2(r+1) + \sum_{r=0}^\infty \frac{(\lambda_\beta t)^{r+1}}{r!} \lambda_\alpha \right]
					\notag \\
				& = \lambda_\alpha \lambda_\beta t + \lambda_\alpha^2 \lambda_\beta t + \lambda_\alpha^2
				\lambda_\beta^2 t^2. \notag
				\end{align}
				Therefore
				\begin{align}
					\mathbb{V}\text{ar} \hat{N}(t) = \mathbb{E} \hat{N}^2(t) - \left[ \mathbb{E}\hat{N}(t) \right]^2
					= \lambda_\alpha \lambda_\beta (1+\lambda_\alpha) t.
				\end{align}
				This result can be confirmed by applying the Wald's formula for the variance.
				\begin{align}
					\mathbb{V}\text{ar} \hat{N}(t) = \mathbb{V}\text{ar} N_\beta(t) \left[ 
					\mathbb{E} X \right]^2 + \mathbb{V}\text{ar} X \mathbb{E} N_\beta(t)
					= \lambda_\alpha \lambda_\beta (1+\lambda_\alpha)t,
				\end{align}
				because $\mathbb{E}X = \mathbb{V}\text{ar} X = \lambda_\alpha$.
			\end{Remark}
		
			\begin{Remark}
				The state probabilities $\hat{p}_k(t) = \text{Pr} \{ \hat{N}(t)=k \}$ satisfy
				the difference-differential equations
				\begin{align}
					\label{elena2}
					\frac{\mathrm d}{\mathrm dt} \hat{p}_k(t) = -\lambda_\beta \hat{p}_k(t) + \lambda_\beta
					e^{-\lambda_\alpha} \sum_{m=0}^k \frac{\lambda_\alpha^m}{m!} \hat{p}_{k-m}(t), \qquad k \geq 0.
				\end{align}
				In order to prove \eqref{elena2} we write
				\begin{align}
					\hat{p}_k(t) = \sum_{r=0}^\infty \text{Pr} \{ N_\alpha(r) = k \} \text{Pr} \{ N_\beta(t) = r \},
				\end{align}
				so that
				\begin{align}
					\frac{\mathrm d}{\mathrm dt} \hat{p}_k(t) & = \sum_{r=0}^\infty \text{Pr} \{ N_\alpha(r) = k \}
					\frac{\mathrm d}{\mathrm dt} \text{Pr} \{ N_\beta(t) = r \} \\
					& = -\lambda_\beta \sum_{r=0}^\infty \text{Pr} \{ N_\alpha(r) = k \}
					\text{Pr} \{ N_\beta(t)=r \} + \lambda_\beta \sum_{r=0}^\infty \text{Pr} \{ N_\alpha(r)=k \}
					\text{Pr} \{ N_\beta(t) = r-1 \} \notag \\
					& = - \lambda_\beta \hat{p}_k(t) + \lambda_\beta \frac{\lambda_\alpha^k}{k!} e^{-\lambda_\beta t}
					\sum_{r=1}^\infty e^{-\lambda_\alpha r} \frac{r^k(\lambda_\beta t)^{r-1}}{(r-1)!} \notag \\
					& = - \lambda_\beta \hat{p}_k(t) + \lambda_\beta \frac{\lambda_\alpha^k}{k!} e^{-\lambda_\beta t}
					\sum_{r=0}^\infty e^{-\lambda_\alpha(r+1)} \frac{(r+1)^k(\lambda_\beta t)^r}{r!} \notag \\
					& = - \lambda_\beta \hat{p}_k(t) + \lambda_\beta e^{-\lambda_\beta t} \frac{\lambda_\alpha^k}{k!}
					e^{-\lambda_\alpha}
					\sum_{m=0}^k \frac{k!}{m!(k-m)!} \sum_{r=0}^\infty e^{-\lambda_\alpha r}
					r^m \frac{(\lambda_\beta t)^r}{r!} \notag \\
					& = - \lambda_\beta \hat{p}_k(t) + \lambda_\beta \lambda_\alpha^k e^{-\lambda_\alpha} \sum_{m=0}^k
					\frac{\lambda_\alpha^{-m}}{(k-m)!} \hat{p}_m(t) \notag \\
					& = - \lambda_\beta \hat{p}_k(t) +
					\lambda_\beta e^{-\lambda_\alpha} \sum_{m=0}^k \frac{\lambda_\alpha^m}{m!} \hat{p}_{k-m}(t). \notag 
				\end{align}
				It is now easy to show that
				\begin{align}
					\hat{G}(u,t) = \sum_{k=0}^\infty u^k \hat{p}_k(t)
				\end{align}
				satisfies the partial differential equation
				\begin{align}
					\frac{\partial}{\partial t} \hat{G}(u,t) = -\lambda_\beta \hat{G}(u,t) + \lambda_\beta
					e^{\lambda_\alpha(u-1)} \hat{G} (u,t) = \lambda_\beta \hat{G} (u,t) \left(
					e^{\lambda_\alpha(u-1)}-1 \right).
				\end{align}
			\end{Remark}

			\begin{Remark}
				For the composition $N_\alpha(\mathfrak{N}(t))$ we have that
				\begin{align}
					\label{ielena}
					\sum_{k=0}^\infty u^k \text{Pr} \{ N_\alpha(\mathfrak{N}(t)) =k \} & =
					\sum_{k=0}^\infty u^k \sum_{r=0}^\infty \frac{(\lambda_\alpha r)^k}{k!}
					e^{-\int_0^k \lambda(w) \mathrm dw} \frac{\left[ \int_0^t \lambda(w) \mathrm dw \right]^r}{r!}
					e^{-\lambda_\alpha r} \\
					& = \sum_{r=0}^\infty e^{-\lambda_\alpha r} \frac{\left[ \int_0^t \lambda(w) \mathrm dw \right]^r}{r!}
					e^{\lambda_\alpha r u} e^{-\int_0^t \lambda(w) \mathrm dw} \notag \\
					& = e^{\left[ e^{\lambda_\alpha(u-1)} -1 \right] \int_0^t \lambda(w) \mathrm dw}. \notag
				\end{align}
				We can also ascertain that
				\begin{align}
					N_\alpha(\mathfrak{N}(t)) \overset{\text{d}}{=} \sum_{j=1}^{\mathfrak{N}(t)}X_j
				\end{align}
				by using result \eqref{ielena}.
			\end{Remark}
			
			\subsubsection{Hitting time for the iterated Poisson process}
			
				Here we study the distribution of the random variable
				\begin{align}
					T_k = \inf (s\colon N_\alpha(N_\beta(s))=k ), \qquad k \geq 1,
				\end{align}
				which represents the first-passage time of the iterated Poisson process at level $k$.
				In the next theorem we state the main result.

				\begin{Theorem}
					\begin{align}
						\label{garla}
						\text{Pr} \{ T_k \in \mathrm ds \} = \mathrm ds \lambda_\beta e^{-\lambda_\alpha}
						\frac{\lambda_\alpha^k}{k!} e^{-\lambda_\beta s} \sum_{j=0}^\infty e^{-\lambda_\alpha j}
						\left[ (j+1)^k - j^k \right] \frac{(\lambda_\beta s)^j}{j!}, \qquad s>0.
					\end{align}
					
					\begin{proof}
						In order to arrive at result \eqref{garla}, we first write
						\begin{align}
							\label{ach}
							\text{Pr} \{ T_k \in \mathrm ds \} & = \sum_{h=1}^k \text{Pr} \{ N_\alpha(N_\beta(s))
							= k-h, N_\alpha(N_\beta(s+\mathrm ds)) = k \} \\
							& = \sum_{h=1}^k \text{Pr} \{ N_\alpha(N_\beta(s)) = k-h, N_\alpha
							(N_\beta(s) + \mathrm d N_\beta(s)) = k \}. \notag
						\end{align}
						Clearly, $\mathrm d N_\beta(s)$ either takes the value $0$ with probability $1-\lambda_\beta
						\mathrm ds$ (and, in this case, all events appearing in \eqref{ach} are mutually
						exclusive) or the value 1. Therefore
						\begin{align}
							\text{Pr} \{ T_k \in \mathrm ds \} = \lambda_\beta \mathrm ds \sum_{h=1}^k
							\text{Pr} \{ N_\alpha(N_\beta(s)) = k-h, N_\alpha(N_\beta(s)+1) = k \}.
						\end{align}
						In the interval $(N_\beta(s), N_\beta(s)+1)$, the external process can take all possible
						values $0 \leq k-h \leq k-1$. The value $k$ must be excluded because if at time
						$s$, $N_\alpha(N_\beta(s))=k$, the first attainment of value $k$ cannot be
						recorded during the interval $(s,s+\mathrm ds]$. Furthermore
						\begin{align}
							\label{freccia}
							& \text{Pr} \{ N_\alpha(N_\beta(s)) = k-h, N_\alpha(N_\beta(s)+1)=k | N_\beta(s)=j \} \\
							& = \text{Pr} \{ N_\alpha(j) = k-h, N_\alpha(j+1)= k \} \notag \\
							& = \text{Pr} \{ N_\alpha(j) = k-h \} \text{Pr} \{ N_\alpha(j+1)-N_\alpha(j)=h \} \notag \\
							& = e^{-\lambda_\alpha j} \frac{(\lambda_\alpha j)^{k-h}}{(k-h)!} e^{-\lambda_\alpha}
							\frac{\lambda_\alpha^h}{h!}. \notag
						\end{align}
						By inserting \eqref{freccia} into \eqref{ach}, we arrive at
						\begin{align}
							\text{Pr} \{ T_k \in \mathrm ds \} & = \mathrm ds \,
							\lambda_\beta e^{-\lambda_\beta s} \sum_{j=0}^\infty \sum_{h=1}^k
							e^{-\lambda_\alpha j} \frac{(\lambda_\alpha j)^{k-h}}{(k-h)!} e^{-\lambda_\alpha}
							\frac{\lambda_\alpha^h}{h!} \frac{(\lambda_\beta s)^j}{j!} \\
							& = \mathrm ds \, \lambda_\beta e^{-\lambda_\alpha} e^{-\lambda_\beta s}
							\frac{\lambda_\alpha^k}{k!} \sum_{j=0}^\infty e^{-\lambda_\alpha j} \left[ (j+1)^k-j^k \right]
							\frac{(\lambda_\beta s)^j}{j!}. \notag
						\end{align}
						We note that $\text{Pr} \{ T_k < \infty \} < 1$, for all $k\geq 1$.
						From \eqref{garla}, we have that
						\begin{align}
							\text{Pr} \{ T_k<\infty \} & = e^{-\lambda_\alpha} \frac{\lambda_\alpha^k}{k!}
							\sum_{j=0}^\infty e^{-\lambda_\alpha j} \left[ (j+1)^k-j^k \right] \\
							& = e^{-\lambda_\alpha} \frac{\lambda_\alpha^k}{(k-1)!} \sum_{j=0}^\infty
							e^{-\lambda_\alpha j} \int_j^{j+1} x^{k-1} \mathrm dx \notag \\
							& = \frac{\lambda_\alpha^k}{(k-1)!} \sum_{j=0}^\infty \int_j^{j+1} e^{-\lambda_\alpha
							(j+1)} x^{k-1} \mathrm dx \notag \\
							& < \frac{\lambda_\alpha^k}{(k-1)!} \int_0^\infty e^{-\lambda_\alpha x} x^{k-1}
							\mathrm dx = 1. \notag
						\end{align}
					\end{proof}
				\end{Theorem}

				\begin{Remark}
					The previous result shows that there is a positive probability of never reaching level $k$
					for the iterated Poisson process. For some cases this probability can be evaluated
					explicitly.
					\begin{align}
						\text{Pr} \{ T_1<\infty \} = \lambda_\alpha e^{-\lambda_\alpha} \sum_{j=0}^\infty
						e^{-\lambda_\alpha j} = \frac{\lambda_\alpha e^{-\lambda_\alpha}}{1-e^{-\lambda_\alpha}}< 1.
					\end{align}
					This is because
					\begin{align}
						0 < 1-e^{-\lambda_\alpha} -\lambda_\alpha e^{-\lambda_\alpha} = 1- \text{Pr}
						\{ N_\alpha(1)=0 \} - \text{Pr} \{ N_\alpha(1)=1 \}.
					\end{align}
					From \eqref{ach} we can evaluate the distribution of $T_1$ as follows.
					\begin{align}
						\text{Pr} \{ T_1 \in \mathrm ds \} = \mathrm ds \, \lambda_\alpha e^{-\lambda_\alpha}
						\lambda_\beta e^{-\lambda_\beta s (1-e^{-\lambda_\alpha})}, \qquad s>0.
					\end{align}
					By similar calculations we have also that
					\begin{align}
						\text{Pr} \{ T_2 \in \mathrm ds \} = \mathrm ds
						\, \lambda_\beta \frac{\lambda_\alpha^2 e^{-\lambda_\alpha}}{2}
						e^{-\lambda_\beta s (1-e^{-\lambda_\alpha})} \left[ 1+2(\lambda_\beta s)e^{-\lambda_\alpha}
						\right].
					\end{align}
					Furthermore
					\begin{align}
						\text{Pr} \{ T_2 < \infty \} = \left[ \text{Pr} \{ T_1< \infty \} \right]^2
						+ \frac{\lambda_\alpha}{2} \text{Pr} \{ T_1<\infty \}.
					\end{align}
				\end{Remark}
				Finally, an alternative form of \eqref{ach} reads
				\begin{align}
					\text{Pr} \{ T_k \in \mathrm ds \} = \lambda_\beta \mathrm ds \sum_{j=0}^\infty
					\int_j^{j+1} e^{-\lambda_\alpha(j+1-x)} \text{Pr} \{ T_k^\alpha \in \mathrm ds \},
				\end{align}
				where $T_k^\alpha = \inf \{ s\colon N_\alpha(s)=k \}$ and $\text{Pr} \{ T_k^\alpha \in \mathrm ds \}$
				is the Erlang distribution for the external Poisson process.
			
		\subsection{Subordination of a homogeneous Poisson process to a fractional pure birth process}
		
			Let $\mathcal{Y}^\nu(t)$, $t>0$, $0<\nu \leq 1$, be a nonlinear fractional pure birth process with rates
			$\lambda_j>0$, $j \in \mathbb{N}$ (representing an extension of the classical nonlinear
			pure birth process). It has been shown in \citet{polbir} that
			\begin{align}
				\label{law}
				\text{Pr} \{ \mathcal{Y}^\nu(t) = k | \mathcal{Y}^\nu(0) = 1 \} =
				\begin{cases}
					\prod_{j=1}^{k-1} \lambda_j \sum_{m=1}^k \frac{E_{\nu,1}(-\lambda_m t^\nu)}{
					\prod_{l=1,l\neq m}^k (\lambda_l-\lambda_m)}, & k>1, \\
					E_{\nu,1} (-\lambda_1 t^\nu), & k=1,
				\end{cases}
			\end{align}
			where $E_{\nu,1}(x)$ is the Mittag--Leffler function.
			For the process $N_\alpha(\mathcal{Y}^\nu(t))$ we have the following result.
			
			\begin{Theorem}
				For the composition $N_\alpha(\mathcal{Y}^\nu(t))$ we have that
				\begin{align}
					\label{sara}
					N_\alpha(\mathcal{Y}^\nu(t)) \overset{\text{d}}{=} \sum_{j=1}^{\mathcal{Y}^\nu(t)}X_j,
				\end{align}
				where the $X_j$, $j \geq 1$, are i.i.d.\ Poisson random variables with rate $\lambda_\alpha$.
				\begin{proof}
					In order to prove \eqref{sara}, we evaluate the probability generating function of both members.
					\begin{align}
						\label{sara2}
						\mathbb{E}u^{N_\alpha(\mathcal{Y}^\nu(t))} & =
						\sum_{k=0}^\infty u^k \left[ \frac{e^{-\lambda_\alpha} \lambda_\alpha^k}{k!}
						E_{\nu,1} (-\lambda_1 t^\nu)
						+ \sum_{r=2}^\infty e^{-\lambda_\alpha r} \frac{(\lambda_\alpha r)^k}{k!}
						\prod_{j=1}^{r-1} \lambda_j \sum_{m=1}^r \frac{1}{\prod_{l=1,l\neq m}^r(\lambda_l-\lambda_m)}
						E_{\nu,1}(-\lambda_m t^\nu) \right] \\
						& = e^{\lambda_\alpha(u-1)} E_{\nu,1}
						(-\lambda_1 t^\nu)
						+ \sum_{r=2}^\infty \left[ e^{\lambda_\alpha(u-1)} \right]^r
						\prod_{j=1}^{r-1} \lambda_j \sum_{m=1}^r \frac{1}{\prod_{l=1,l\neq m}^r(\lambda_l-\lambda_m)}
						E_{\nu,1}(-\lambda_m t^\nu) . \notag
					\end{align}
					The probability generating function of the right-hand side of \eqref{sara} reads
					\begin{align}
						\label{sara3}
						\mathbb{E} u^{\sum_{j=1}^{\mathcal{Y}^\nu(t)}X_j} =
						e^{\lambda_\alpha(u-1)} \text{Pr} \{ \mathcal{Y}^\nu(t) = 1 \} +
						\sum_{r=2}^\infty e^{\lambda_\alpha(u-1)r} \text{Pr} \{ \mathcal{Y}^\nu(t) = r \},
					\end{align}
					as $\mathbb{E}u^X = e^{\lambda_\alpha(u-1)}$, because $X$ is a Poisson random
					variable of parameter $\lambda_\alpha$.
					By comparing \eqref{sara2} with \eqref{sara3}, in view of \eqref{law}, we have the claimed result.
				\end{proof}
			\end{Theorem}
			
			\begin{Remark}
				For the fractional linear birth process $Y^\nu(t)$, $t>0$, the distribution \eqref{law}
				specialises and takes the form
				\begin{align}
					\label{lawlin}
					\text{Pr} \{ Y^\nu(t) = k | Y^\nu(0) = 1 \} =
					\sum_{j=1}^k \binom{k-1}{j-1} (-1)^{j-1} E_{\nu,1}(-\lambda_\beta j t^\nu).
				\end{align}
				For $\nu=1$, the distribution \eqref{lawlin} reduces to the geometric distribution
				of the classical Yule--Furry process. The probability generating function of the
				composed process $N_\alpha(Y^\nu(t))$ becomes
				\begin{align}
					\mathbb{E} \left[ u^{N_\alpha(Y^\nu(t))} \right] & =
					\sum_{k=0}^\infty u^k \sum_{r=1}^\infty e^{-\lambda_\alpha r} \frac{(\lambda_\alpha r)^k}{k!}
					\sum_{j=1}^r \binom{r-1}{j-1} (-1)^{j-1} E_{\nu,1}(-\lambda_\beta j t^\nu) \\
					& = \sum_{r=1}^\infty e^{-\lambda_\alpha r} e^{\lambda_\alpha r u} \sum_{j=1}^r
					\binom{r-1}{j-1} (-1)^{j-1} E_{\nu,1}(-\lambda_\beta j t^\nu) \notag \\
					& = \sum_{j=1}^\infty (-1)^{j-1} E_{\nu,1} (-\lambda_\beta j t^\nu) \sum_{r=j}^\infty
					\binom{r-1}{j-1} e^{-\lambda_\alpha(1-u) r} \notag \\
					& = \sum_{j=1}^\infty (-1)^{j-1} E_{\nu,1} (-\lambda_\beta j t^\nu) e^{-\lambda_\alpha (i-u)j}
					\sum_{r=0}^\infty \binom{j+r-1}{r} e^{-\lambda_\alpha(1-u)r} \notag \\
					& = \sum_{j=1}^\infty (-1)^{j-1} E_{\nu,1} (-\lambda_\beta j t^\nu) e^{-\lambda_\alpha (1-u)j}
					\sum_{r=0}^\infty \binom{-j}{r} (-1)^r e^{-\lambda_\alpha (1-u)r} \notag \\
					& = \sum_{j=1}^\infty (-1)^{j-1} e^{-\lambda_\alpha(1-u)j} \left( 1-e^{-\lambda_\alpha(1-u)}
					\right)^{-j}
					E_{\nu,1}(-\lambda_\beta j t^\nu) \notag \\
					& = \sum_{j=1}^\infty (-1)^{j-1} \left[ \frac{e^{\lambda_\alpha(u-1)}}{
					1-e^{\lambda_\alpha(u-1)}} \right]^j
					E_{\nu,1}(-\lambda_\beta j t^\nu). \notag
				\end{align}
			\end{Remark}
			
		\subsection{Subordination with a Poisson field}
			
			Let $B$ a Borel set and let us indicate with $\Lambda(B)$ its Lebesgue measure.
			The aim of this section is to study the process
			\begin{align}
				N_\alpha(\mathcal{F}(B)),
			\end{align}
			where $N_\alpha(t)$, $t>0$, is a homogeneous Poisson process with rate $\lambda_\alpha>0$, and
			$\mathcal{F}(B)$, is a homogeneous Poisson field with rate $\lambda>0$.
			Analogously to the iterated Poisson process, the representation
			\begin{align}
				N_\alpha(\mathcal{F}(B)) \overset{\text{d}}{=} X_1+\dots+X_{\mathcal{F}(B)},
			\end{align}
			where the random variables $X_j$ are i.i.d. Poisson distributed of parameter $\lambda_\alpha>0$,
			is still valid. Therefore, the state probabilities can be written as
			\begin{align}
				\text{Pr} \{ N_\alpha(\mathcal{F}(B)) = k \} =
				\frac{\lambda_\alpha^k}{k!} e^{-\lambda \Lambda(B)(1-e^{-\lambda_\alpha})}
				\mathfrak{B}_k \left( \lambda \Lambda(B e^{-\lambda_\alpha}) \right),
				\quad k \geq 0,
			\end{align}
			and the probability generating function reads
			\begin{align}
				G(u,B) = e^{\lambda \Lambda(B)\left( e^{\lambda_\alpha
				\left(u-1 \right)} -1 \right)}, \qquad |u|\leq 1.
			\end{align}
			The emptiness probability is given by:
			\begin{align}
				\text{Pr} \{ N_\alpha(\mathcal{F}(B)) = 0 \} =
				e^{-\lambda\Lambda(B)\left( 1-e^{-\lambda_\alpha} \right)}.
			\end{align}
			Let now $B_l$ be the disc with radius $l$, centred in the origin. The first-contact
			distribution $H(l)$, $l \in \mathbb{R}^+$, (with the first point) is
			\begin{align}
				H(l) & = 1 - \text{Pr} \{ N_\alpha(\mathcal{F}(B_l)) = 1 \}
				= 1 - e^{-\lambda \pi l^2 (1-e^{-\lambda_\alpha})}, \qquad l \in \mathbb{R}^+.
			\end{align}
			The probability density is in turn
			\begin{align}
				h(l) = 2\lambda \pi l (1-e^{-\lambda_\alpha}) e^{-\lambda \pi l^2 (1-e^{-\lambda_\alpha})},
				\qquad l \in \mathbb{R}^+,
			\end{align}
			that is, a Rayleigh distribution (see, for the classical non subordinated case,
			\citet{stoyan}, page 213).
			
	\section{Compositions of Poisson processes with first-passage time of different point processes}
	
		In this section, we consider a fractional Poisson process $N^\nu(t)$ (see \citet{beg} for information on this
		process) whose first-passage time
		\begin{align}
			\label{cusc}
			\tau^\nu_k = \inf (t \colon N^\nu(t) = k), \qquad k \geq 1,
		\end{align}
		is composed either with a homogeneous Poisson process $N_\alpha(t)$ ($\lambda_\alpha>0$ is its rate) or
		with a Yule--Furry process $Y_\alpha(t)$.
		
		The distribution of \eqref{cusc} has the following density
		\begin{align}
			\text{Pr}(\tau^\nu_k \in \mathrm dt)/\mathrm dt & = \frac{\mathrm d}{\mathrm dt}
			\text{Pr} \{ N^\nu(t) \geq k \}
			= \lambda_\beta^k t^{\nu k-1}
			E_{\nu,\nu k}^k (-\lambda_\beta t^\nu), \qquad t>0, \: \lambda_\beta > 0,
		\end{align}
		\citep[formula (1.6)]{beg2}, where
		\begin{align}
			E_{\xi,\gamma}^\delta(z) = \sum_{r=0}^\infty \frac{(\delta)_r z^r}{\Gamma(\xi r +
			\gamma) r!}, \qquad \xi, \gamma, \delta \in \mathbb{C}, \: \mathfrak{R}(\xi)>0.
		\end{align}
		The function $E_{\xi,\gamma}^\delta(z)$ is called generalised Mittag--Leffler function
		(see \citet[page 91]{haubold})
		
		\begin{Theorem}
			The composed process $\tilde{N}^\nu(k) = N_\alpha(\tau_k^\nu)$ has the distribution
			\begin{align}
				\label{car2new}
				\text{Pr} \{ \tilde{N}^\nu(k) = r \}
				= \frac{1}{r!} \sum_{j=0}^\infty \binom{-k}{j} \left( \frac{\lambda_\beta}{
				\lambda_\alpha^\nu} \right)^{k+j} \frac{\Gamma(\nu(k+j)+r)}{\Gamma(\nu(k+j))}, \qquad r\geq 0,
			\end{align}
			and possesses probability generating function
			\begin{align}
				\label{car3new}
				\mathbb{E}u^{\tilde{N}^\nu(k)}
				= \left[ 1+(1-u)^\nu \frac{\lambda_\alpha^\nu}{\lambda_\beta} \right]^{-k}, \qquad |u|\leq 1.
			\end{align}
			\begin{proof}
				We start our proof with the following relation:
				\begin{align}
					\label{caratterenew}
					\text{Pr} \left\{ \tilde{N}^\nu(k) = r \right\} = \int_0^\infty \text{Pr} \left\{ N_\alpha(s)
					=r \right\} \text{Pr} \left\{ \tau_k^\nu \in \mathrm ds \right\}, \qquad k \geq 1, \: r \geq 0.
				\end{align}
				Instead of writing the distribution of $\tau_k^\nu$ in terms of generalised Mittag--Leffler
				functions, it is more convenient to work with the following expression which we derive from scratch.
				
				From \citet{beg}, we know that the distribution of the fractional Poisson process is
				\begin{align}
					\text{Pr} \{ N^\nu(t) = m \} & = \sum_{j=m}^\infty (-1)^{j-m} \binom{j}{m} \frac{(\lambda_\beta
					t^\nu)^j}{\Gamma(\nu j +1)} \\
					& = (\lambda_\beta t^\nu)^m
					\sum_{j=0}^\infty \binom{j+m}{j} \frac{(-\lambda_\beta t^\nu)^j}{\Gamma \left(\nu(m+j)+1\right)},
					\qquad t>0, \: \nu \in (0,1]. \notag
				\end{align}
				By considering that $\text{Pr} \{ \tau_k^\nu < s \} = \text{Pr} \{ N^\nu(s) \geq k \}$, we have
				that
				\begin{align}
					\label{car1new}
					\text{Pr} \{ \tau_k^\nu \in \mathrm ds \}/\mathrm ds & = \frac{\mathrm d}{\mathrm ds} \sum_{m=k}^\infty
					\text{Pr} \{ N^\nu(s) = m \} \\
					& = \frac{\mathrm d}{\mathrm ds} \sum_{m=k}^\infty (\lambda_\beta s^\nu)^m \sum_{j=0}^\infty
					\binom{j+m}{j} \frac{(-\lambda_\beta t^\nu)^j}{\Gamma \left(\nu(m+j)+1\right)} \notag \\
					& = \frac{\mathrm d}{\mathrm ds} \sum_{m=k}^\infty \lambda_\beta^m \sum_{j=0}^\infty
					\binom{j+m}{j} (-\lambda_\beta)^j \frac{(s^\nu)^{m+j}}{\Gamma(\nu(m+j)+1)} \notag \\
					& \overset{h=j+m}{=} \frac{\mathrm d}{\mathrm ds} \sum_{m=k}^\infty \lambda_\beta^m
					\sum_{h=m}^\infty \binom{h}{h-m} (-\lambda_\beta)^{h-m} \frac{s^{\nu h}}{\Gamma(\nu h +1)} \notag \\
					& = \frac{\mathrm d}{\mathrm ds} \sum_{h=k}^\infty (-\lambda_\beta)^h \frac{s^{\nu h}}{
					\Gamma(\nu h +1)} \sum_{m=k}^h \binom{h}{m} (-1)^m \notag \\
					& = \frac{\mathrm d}{\mathrm ds} \sum_{h=k}^\infty (-\lambda_\beta)^h \frac{s^{\nu h}}{
					\Gamma(\nu h+1)} (-1)^k \binom{h-1}{k-1} \notag \\
					& = \sum_{h=k}^\infty (-\lambda_\beta)^h (-1)^k \binom{h-1}{k-1} \frac{\nu h s^{\nu h -1}}{
					\Gamma(\nu h +1)} \notag \\
					& = \sum_{h=k}^\infty (-\lambda_\beta)^h (-1)^k \binom{h-1}{k-1} \frac{s^{\nu h -1}}{\Gamma(\nu h)}
					\notag \\
					& \overset{j=h-k}{=} \sum_{j=0}^\infty (-\lambda_\beta)^{j+k} (-1)^k \binom{j+k-1}{k-1}
					\frac{s^{\nu(k+j)-1}}{\Gamma(\nu(k+j))} \notag \\
					& = \lambda_\beta^k \sum_{j=0}^\infty \binom{-k}{j} \lambda_\beta^{j} \frac{s^{\nu(k+j)-1}}{
					\Gamma(\nu(k+j))}. \notag
				\end{align}
				In the fifth step of \eqref{car1new}, we used the following formula which is
				interesting in itself
				\begin{align}
					\label{starnew}
					\sum_{m=k}^h \binom{h}{m} (-1)^m = (-1)^k \binom{h-1}{k-1}.
				\end{align}
				We provide here a proof of \eqref{starnew} in the following way:
				\begin{align}
					\sum_{m=k}^h & \binom{h}{m} (-1)^m
					= \sum_{m=0}^h \binom{h}{m} (-1)^m - \sum_{m=0}^{k-1}
					\binom{h}{m} (-1)^m
					= \sum_{m=0}^{k-1} \binom{h}{m}(-1)^{m+1} \\
					= {} & -1+h-\frac{h(h-1)}{2}+\frac{h(h-1)(h-2)}{2\cdot 3} - \frac{h(h-1)(h-2)(h-3)}{2\cdot 3 \cdot 4}
					\notag \\
					& + \dots + (-1)^k \frac{h(h-1)(h-2)\dots (h-k+2)}{2 \cdot 3 \cdots (k-1)} \notag \\
					= {} & (h-1) \left[ 1- \frac{h}{2} + \frac{h(h-2)}{2\cdot 3} - \frac{h(h-2)(h-3)}{2\cdot 3\cdot 4}
					\right. \notag \\
					& \left. +\dots + (-1)^k \frac{h(h-2)\dots (h-k+2)}{2\cdot 3 \cdots (k-1)} \right] \notag \\
					= {} & \frac{(h-1)(h-2)}{2} \left[ -1 + \frac{h}{3} - \frac{h(h-3)}{3\cdot 4}
					+\dots +(-1)^k \frac{h(h-3)\dots (h-k+2)}{3 \cdot 4 \cdots (k-1)} \right] \notag \\
					= {} & \frac{(h-1)(h-2)(h-3)}{2 \cdot 3} \left[ 1-\frac{h}{4} +
					\dots + (-1)^k \frac{h(h-4)\dots (h-k+2)}{4\cdot 5\cdots (k-1)}\right] = \dots = \notag \\
					= {} & \frac{(h-1)(h-2)(h-3)\dots(h-k+2)}{(k-2)!} \left[ (-1)^{k-1} + (-1)^k \frac{h}{k-1} \right]
					\notag \\
					= {} & \frac{(h-1)(h-2)(h-3)\dots(h-k+2)}{(k-2)!} (-1)^k \left[ -1+\frac{h}{k-1} \right] \notag \\
					= {} & (-1)^k \frac{(h-1)!}{(k-1)!(h-k)!}
					= (-1)^k \binom{h-1}{k-1}. \notag
				\end{align}
				By inserting \eqref{car1new} into \eqref{caratterenew} we obtain that
				\begin{align}
					\text{Pr} \{ \tilde{N}^\nu(k) = r \}
					& = \int_0^\infty e^{-\lambda_\alpha s} \frac{(\lambda_\alpha s)^r}{r!}
					\lambda_\beta^k \sum_{j=0}^\infty \binom{-k}{j} \lambda_\beta^j
					\frac{s^{\nu(k+j)-1}}{\Gamma(\nu(k+j))} \mathrm ds \\
					& = \frac{\lambda_\alpha^r}{r!} \lambda_\beta^k \sum_{j=0}^\infty \binom{-k}{j} \lambda_\beta^j
					\frac{\Gamma(\nu(k+j)+r)}{\lambda_\alpha^{r+\nu(k+j)}\Gamma(\nu(k+j))} \notag \\
					& = \frac{1}{r!} \sum_{j=0}^\infty \binom{-k}{j} \left( \frac{\lambda_\beta}{
					\lambda_\alpha^\nu} \right)^{k+j} \frac{\Gamma(\nu(k+j)+r)}{\Gamma(\nu(k+j))}. \notag
				\end{align}
				This proves formula \eqref{car2new}.
				The probability generating function
				\begin{align}
					\mathbb{E}u^{\tilde{N}^\nu(k)}
					& = \sum_{r=0}^\infty u^r \int_0^\infty \frac{e^{-\lambda_\alpha s}(\lambda_\alpha s)^r}{r!}
					\lambda_\beta^k \sum_{j=0}^\infty \binom{-k}{j} \lambda_\beta^j
					\frac{s^{\nu(k+j)-1}}{\Gamma(\nu(k+j))} \mathrm ds \\
					& = \int_0^\infty e^{-\lambda_\alpha s} e^{\lambda_\alpha s u} \lambda_\beta^k \sum_{j=0}^\infty
					\binom{-k}{j} \lambda_\beta^j	\frac{s^{\nu(k+j)-1}}{\Gamma(\nu(k+j))} \mathrm ds \notag \\
					& =\left( \frac{\lambda_\beta}{\lambda_\alpha^\nu(1-u)^\nu} \right)^k
					\sum_{j=0}^\infty \binom{-k}{j} \left( \frac{\lambda_\beta}{\lambda_\alpha^\nu(1-u)^\nu} \right)^j
					\notag \\
					& = \left( \frac{\lambda_\beta}{\lambda_\alpha^\nu(1-u)^\nu} \right)^k \left(
					1+ \frac{\lambda_\beta}{\lambda_\alpha^\nu(1-u)^\nu} \right)^{-k} \notag \\
					& = \left[ 1+(1-u)^\nu \frac{\lambda_\alpha^\nu}{\lambda_\beta} \right]^{-k}. \notag
				\end{align}
			\end{proof}
		\end{Theorem}
		
		\begin{Remark}
			The waiting time $\tau_k^\nu$, $k\geq 1$, for the kth event of the fractional Poisson process
			can be viewed as the sum of independent waiting times $\tau_{1,j}^\nu$ separating the
			events of the Poisson flow, i.e.
			\begin{align}
				\label{prisma}
				\tau_k^\nu = \sum_{j=1}^k \tau_{1,j}^\nu.
			\end{align}			
			It is well-known (see e.g.\ \citet{mainardi} or \citet{beg})
			\begin{align}
				\label{unaparola}
				\text{Pr} \{ \tau_{1,j}^\nu \in \mathrm ds \} & = \text{Pr} \{ \tau_1^\nu \in \mathrm ds \} \\
				& = \lambda_\beta s^{\nu-1} \sum_{j=0}^\infty \frac{(-\lambda_\beta s^\nu)^j}{\Gamma(\nu(j+1))} \notag \\
				& = \lambda_\beta s^{\nu-1} E_{\nu,\nu} (-\lambda_\beta s^\nu) \notag \\
				& = - \frac{\mathrm d}{\mathrm ds} E_{\nu,1} (-\lambda_\beta s^\nu). \notag
			\end{align}
			From \eqref{unaparola} it follows that
			\begin{align}
				\int_0^\infty \text{Pr} \{ \tau_{1,j}^\nu \in \mathrm ds \} = 1,
			\end{align}
			and, by writing the distribution of $\tau_k^\nu$ as convolution of the terms pertaining to $\tau_{1,j}^\nu$,
			we have also that
			\begin{align}
				\int_0^\infty \text{Pr} \{ \tau_k^\nu \in \mathrm ds \} = 1.
			\end{align}
			The Laplace transform of \eqref{car1new} is easily calculated and reads
			\begin{align}
				\label{ragno}
				\mathbb{E} e^{-\mu \tau_k^\nu} & = \int_0^\infty e^{-\mu s} 
				\lambda_\beta^k \sum_{j=0}^\infty \binom{-k}{j} \lambda_\beta^{j} \frac{s^{\nu(k+j)-1}}{
				\Gamma(\nu(k+j))} \mathrm ds \\
				& = \frac{\lambda_\beta^k}{\mu^{\nu k}} \sum_{j=0}^\infty \binom{-k}{j}
				\left( \frac{\lambda_\beta}{\mu^\nu}\right)^j 
				= \left[ \frac{\mu^\nu}{\lambda_\beta} +1 \right]^{-k}, \notag
			\end{align}			
			and this clearly confirms the additive structure \eqref{prisma}.
		\end{Remark}

		\begin{Remark}
			Result \eqref{ragno} suggests the following representation for the composed process
			$\tilde{N}^\nu(k)$:
			\begin{align}
				\tilde{N}^\nu(k) \overset{\text{d}}{=} \xi_1^\nu + \dots + \xi_k^\nu,
			\end{align}
			where the $\xi_j^\nu$ are independent random variables which are called discrete Mittag--Leffler
			random variables (see \citet{pillai}) having parameters $\nu$ and $\lambda_\alpha^\nu/\lambda_\beta$.
			The discrete Mittag--Leffler reduces to the geometric random variable of parameter
			$q=\lambda_\alpha/(\lambda_\alpha+\lambda_\beta)$ when $\nu=1$. We now give the explicit distribution
			of the generalised geometric random variables $\xi_j^\nu$.
			\begin{align}
				\label{dmlnew}
				\text{Pr} \{ \xi^\nu =r \} & = \frac{1}{r!} \sum_{j=0}^\infty
				\left( \frac{\lambda_\beta}{\lambda_\alpha^\nu} \right)^{j+1}
				\frac{\Gamma(\nu(j+1)+r)}{\Gamma(\nu(j+1))} \\
				& = \frac{\lambda_\beta}{\lambda_\alpha^\nu r!} \int_0^\infty e^{-w}
				\sum_{j=0}^\infty \frac{w^{\nu(j+1)+r-1}}{\Gamma(\nu(j+1))}
				\left( - \frac{\lambda_\beta}{\lambda_\alpha^\nu} \right)^j \mathrm dw \notag \\
				& = \frac{\lambda_\beta}{\lambda_\alpha^\nu} \int_0^\infty e^{-w}
				\frac{w^{\nu+r-1}}{r!} E_{\nu,\nu} \left( -\frac{\lambda_\beta}{\lambda_\alpha^\nu}
				w^\nu \right) \mathrm dw. \notag
			\end{align}
			For $\nu=1$, we extract from \eqref{dmlnew} the geometric distribution.
			\begin{align}
				\text{Pr} \{ \xi^1 = r \} & = \frac{\lambda_\beta}{\lambda_\alpha} \int_0^\infty
				e^{-w} \frac{w^r}{r!} e^{-\frac{\lambda_\beta w}{\lambda_\alpha}w} \mathrm dw
				= \frac{\lambda_\beta/\lambda_\alpha}{\left(1+\frac{\lambda_\beta}{\lambda_\alpha}\right)^r}
				= pq^{r-1},
			\end{align}
			where $p=\lambda_\beta/(\lambda_\alpha+\lambda_\beta)$.
			In order to check that the generalised geometric law \eqref{dmlnew} sums up to unity, we write
			\begin{align}
				\sum_{r=0}^\infty \text{Pr} \{ \xi^\nu = r \} & = \frac{\lambda_\beta}{\lambda_\alpha^\nu}
				\int_0^\infty e^{-w} w^{\nu-1} E_{\nu,\nu} \left( -\frac{\lambda_\beta}{\lambda_\alpha^\nu}
				\right) \mathrm dw \\
				& = \int_0^\infty - \frac{\mathrm d}{\mathrm dw} E_{\nu,1} \left( -\frac{\lambda_\beta}{\lambda_\alpha^\nu}
				\right) \mathrm dw
				= \left| -E_{\nu,1} \left( -\frac{\lambda_\beta}{\lambda_\alpha^\nu}
				\right) \right|_{w=0}^{w=\infty} = 1. \notag
			\end{align}
			
			Furthermore, formula \eqref{car3new} shows that $\tilde{N}^\nu(t)$ possesses Linnik distribution
			and its form is explicitly given by \eqref{car2new}.
			
			Finally, the distribution \eqref{car2new}, for $\nu=1$, becomes the negative binomial distribution
			having parameters $k$ and $\lambda_\alpha/(\lambda_\alpha+\lambda_\beta)$. Indeed, from
			\eqref{car2new}, we have that
			\begin{align}
				\label{viewnew}
				\text{Pr} \left(\tilde{N}^1(k) = r\right)
				& = \sum_{m=0}^\infty (-1)^m \left( \frac{\lambda_\beta}{\lambda_\alpha} \right)^{k+m}
				\frac{(k+m+r-1)!}{r!m!(k-1)!} \\
				& = \binom{k+r-1}{r}
				\sum_{m=0}^\infty \left( \frac{\lambda_\beta}{\lambda_\alpha} \right)^{k+m} (-1)^m
				\binom{k+m+r-1}{m} \notag \\
				& = \binom{k+r-1}{r} \left( \frac{\lambda_\beta}{\lambda_\alpha} \right)^k
				\sum_{m=0}^\infty \left( \frac{\lambda_\beta}{\lambda_\alpha} \right)^m
				\binom{-(k+r)}{m} \notag \\
				& = \binom{k+r-1}{r} \left( \frac{\lambda_\beta}{\lambda_\alpha} \right)^k
				\left( 1+\frac{\lambda_\beta}{\lambda_\alpha} \right)^{-\left(k+r\right)} \notag \\
				& = \binom{k+r-1}{r} \left( \frac{\lambda_\alpha}{\lambda_\alpha+\lambda_\beta} \right)^r
				\left( 1-\frac{\lambda_\alpha}{\lambda_\alpha+\lambda_\beta} \right)^k. \notag
			\end{align}
		\end{Remark}
		
		\begin{Remark}
			We now find the distribution of a slightly modified first-passage time
			\begin{align}
				\label{piripo}
				\hat{\tau}_k^\nu = \left( \frac{t}{k} \right)^{1/\nu} \tau_k^\nu,
			\end{align}
			where $\tau_k^\nu$ is defined in \eqref{cusc} and has distribution \eqref{car1new}.
			\begin{align}
				\label{bellabel}
				\text{Pr} \{ \hat{\tau}_k^\nu \in \mathrm ds \} & = \mathrm ds \left( \frac{k}{t} \right)^{1/\nu}
				\lambda_\beta^k \sum_{j=0}^\infty \binom{-k}{j}
				\lambda_\beta^j \frac{\left[ s \left( \frac{k}{t} \right)^{1/\nu}
				\right]^{\nu(k+j)-1}}{\Gamma(\nu(k+j))} \\
				& = \mathrm ds \sum_{j=0}^\infty \binom{-k}{j} \lambda_\beta^{k+j} \left( \frac{k}{t} \right)^{k+j}
				\frac{s^{\nu(k+j)-1}}{\Gamma(\nu(k+j))}. \notag
			\end{align}
			The Laplace transform of \eqref{bellabel} becomes
			\begin{align}
				\int_0^\infty e^{-\mu s} \text{Pr} \{ \hat{\tau}_k^\nu \in \mathrm ds \}
				= \left[ 1 + \lambda_\beta \frac{t \mu^\nu}{k} \right]^{-k}.
			\end{align}
			For $k \rightarrow \infty$ we obtain the fine result
			\begin{align}
				\label{piripi}
				\lim_{k \rightarrow \infty} \int_0^\infty e^{-\mu s} \text{Pr} \{ \hat{\tau}_k^\nu \in \mathrm ds \}
				= e^{-\lambda_\beta t \mu^\nu}.
			\end{align}
			Result \eqref{piripi} shows that the rescaled first-passage time \eqref{piripo} converges in
			distribution to a positively skewed stable law of order $\nu \in (0,1)$.
		\end{Remark}
		
		We now consider the Yule--Furry process $Y_\alpha(t)$, with a single progenitor, subordinated to
		the first-passage time $\tau_k^\nu$.
		The distribution of $Y_\alpha(\tau_k^\nu)$ is given below and can be determined as follows.
		Bearing in mind the distribution \eqref{car1new}, we have that
		\begin{align}
			\label{bottnew}
			\text{Pr} \{ Y_\alpha (\tau_k^\nu)=r \}
			& = \int_0^\infty e^{-\lambda_\alpha s}(1-e^{-\lambda_\alpha s})^{r-1} \lambda_\beta^k
			\sum_{j=0}^\infty \binom{-k}{j} \lambda_\beta^{j} \frac{s^{\nu(k+j)-1}}{
			\Gamma(\nu(k+j))} \mathrm ds \\
			& = \int_0^\infty \sum_{h=1}^r \binom{r-1}{h-1} (-1)^{h-1} e^{-\lambda_\alpha h s}
			\lambda_\beta^k \sum_{j=0}^\infty \binom{-k}{j} \lambda_\beta^{j} \frac{s^{\nu(k+j)-1}}{
			\Gamma(\nu(k+j))} \mathrm ds \notag \\
			& = \lambda_\beta^k \sum_{h=1}^r \binom{r-1}{h-1} (-1)^{h-1} \sum_{j=0}^\infty
			\binom{-k}{j} \lambda_\beta^j \frac{1}{(\lambda_\alpha h)^{\nu(k+j)}} \notag \\
			& \left( \text{if} \: \, \lambda_\beta/\lambda_\alpha^\nu < 1 \right) \notag \\
			& = \sum_{h=1}^r \left[ \frac{\lambda_\beta}{\lambda_\alpha^\nu h^\nu} \right]^k
			\binom{r-1}{h-1} (-1)^{h-1} \left[1+\frac{\lambda_\beta}{\lambda_\alpha^\nu h^\nu} \right]^{-k}
			\notag \\
			& = \sum_{h=1}^r \binom{r-1}{h-1} (-1)^{h-1} \left[ 1+h^\nu\frac{\lambda_\alpha^\nu}{
			\lambda_\beta} \right]^{-k}, \qquad r \geq 0. \notag
		\end{align}			
		The probability generating function of the distribution \eqref{bottnew} reads
		\begin{align}
			\label{prpr}
			\mathbb{E}u^{Y_\alpha (\tau_k^\nu)} & = \sum_{r=1}^\infty u^r \sum_{h=1}^r
			\binom{r-1}{h-1} (-1)^{h-1} \left[ 1+h^\nu \frac{\lambda_\alpha^\nu}{\lambda_\beta} \right]^{-k} \\
			& = \sum_{h=1}^\infty (-1)^{h-1} \left[ 1+h^\nu \frac{\lambda_\alpha^\nu}{\lambda_\beta} \right]^{-k}
			\sum_{r=h}^\infty u^r \binom{r-1}{h-1} \notag \\
			& = \sum_{h=1}^\infty (-1)^{h-1} \left( \frac{u}{1-u} \right)^h \left[ 1+h^\nu
			\frac{\lambda_\alpha^\nu}{\lambda_\beta} \right]^{-k}, \qquad |u|<1. \notag
		\end{align}
		We remark that the inversion of sums in \eqref{prpr} is valid only for $|u|<1$.

		\subsubsection{The classical case $\bm{\nu=1}$}
			
			For $\nu=1$ we have special interesting results for $N_\alpha(\tau_k^1)$ and $Y_\alpha(\tau_k^1)$.
			For the first process we have the following result.
			
			\begin{Theorem}
				The composed process $N_\alpha(\tau_k^1)$ has the following representation:
				\begin{align}
					\label{dueparole}
					N_\alpha(\tau_k^1) \overset{\text{d}}{=} X_1+\dots +X_N, 
				\end{align}
				where $N$ is a Poisson random variable of parameter
				\begin{align}
					\mu=\log \left( \frac{\lambda_\alpha+\lambda_\beta}{\lambda_\beta} \right)^k,
				\end{align}
				and the $X_j$s are i.i.d.\ random variables with logarithmic distribution of parameter
				$q=\lambda_\alpha/(\lambda_\alpha+\lambda_\beta)$.
				\begin{proof}
					The random variable $N_\alpha(\tau_k^1)$ is a negative binomial $W$ (see \eqref{viewnew})
					with distribution
					\begin{align}
						\text{Pr} \{ W = r\} = \binom{k+r-1}{r} p^k q^r.
					\end{align}
					In our case $p=\lambda_\beta/(\lambda_\alpha+\lambda_\beta)$ and
					$q=\lambda_\alpha/(\lambda_\alpha+\lambda_\beta)$.
					It is well-known that it can be expanded as a random sum of the form
					\begin{align}
						N_\alpha(\tau_k^1) \overset{\text{d}}{=} X_1+\dots +X_N, 
					\end{align}
					where $N$ is a Poisson random variable of parameter $\mu=-k\log p$ and
					$X$ is a logarithmic distribution of parameter $q$.
				\end{proof}
			\end{Theorem}
			
			\begin{Remark}
				From \eqref{viewnew} we can infer that
				\begin{align}
					\label{media}
					\mathbb{E} \tilde{N}(k) = \frac{\lambda_\alpha +\lambda_\beta}{\lambda_\beta} k,
				\end{align}
				and
				
				\begin{align}
					\label{varianza}
					\mathbb{V}\text{ar} \tilde{N}(k) = \frac{\lambda_\alpha(\lambda_\alpha+\lambda_\beta)}{
					\lambda_\beta^2} k.
				\end{align}
			\end{Remark}
			
			These results can be confirmed by applying Wald's formula to the random sum \eqref{dueparole}.
				
			\begin{Remark}
				For $\nu=1$, the distribution \eqref{bottnew} becomes
				\begin{align}
					\label{geodetica}
					\text{Pr} \{ Y_\alpha(\tau_k^1) = r \} =
					\sum_{h=1}^r \binom{r-1}{h-1} (-1)^{h-1} \left[ 1+h\frac{\lambda_\alpha}{
					\lambda_\beta} \right]^{-k}, \qquad r \geq 1, \: k \geq 1.
				\end{align}
				We are able to give a fine expression of \eqref{geodetica} for $k=1$. We have
				\begin{align}
					\label{curvatura}				
					\text{Pr} \{ Y_\alpha(\tau_1^1) = r \} & =
					\sum_{h=1}^r \binom{r-1}{h-1} (-1)^{h-1} \frac{\lambda_\beta}{\lambda_\beta+h \lambda_\alpha} \\
					& = \sum_{h=0}^{r-1} \binom{r-1}{h} (-1)^h \frac{\lambda_\beta}{\lambda_\alpha}
					\frac{1}{\left( 1+\frac{\lambda_\beta}{\lambda_\alpha} +h \right)}. \notag
				\end{align}
				In light of the formula	(see e.g.\ \citet{kir})
				\begin{align}
					\label{kir}
					\sum_{k=0}^N \binom{N}{k} (-1)^k \frac{1}{x+k} = \frac{N!}{x(x+1)\cdots(x+N)},
				\end{align}
				the probability \eqref{curvatura} becomes
				\begin{align}
					\label{maus}
					\text{Pr} \{ Y_a(\tau_1^1) =r \} & = \frac{\lambda_\beta}{\lambda_\alpha}
					\frac{(r-1)!\Gamma\left( \frac{\lambda_\beta}{\lambda_\alpha}+1
					\right)}{\Gamma\left( \frac{\lambda_\beta}{\lambda_\alpha}+1+r \right)}
					= \frac{\lambda_\beta}{\lambda_\alpha}
					\text{Beta}\left( r,\frac{\lambda_\beta}{\lambda_\alpha}+1 \right), \qquad r \geq 1.
				\end{align}
				We can easily check that \eqref{maus} sums up to unity because
				\begin{align}
					\sum_{r=1}^\infty \frac{\lambda_\beta}{\lambda_\alpha}
					\text{Beta}\left( r,\frac{\lambda_\beta}{\lambda_\alpha}+1 \right)
					& = \sum_{r=1}^\infty \frac{\lambda_\beta}{\lambda_\alpha}
					\int_0^1 x^{r-1}(1-x)^{\lambda_\beta/\lambda_\alpha} \mathrm dx \\
					& = \frac{\lambda_\beta}{\lambda_\alpha} \int_0^1 (1-x)^{\lambda_\beta/\lambda_\alpha-1}
					=1. \notag
				\end{align}
			\end{Remark}
			
			\begin{Remark}
				The mean value and the variance of $Y_\alpha(\tau_k^1)$ can be obtained by means of the following
				calculations.
				\begin{align}
					\mathbb{E}Y_\alpha(\tau_k^1) & = \frac{\lambda_\beta^k}{(k-1)!} \int_0^\infty e^{\lambda_\alpha s}
					s^{k-1} e^{-\lambda_\beta s} \mathrm ds
					= \frac{\lambda_\beta^k}{(k-1)!} \int_0^\infty e^{-(\lambda_\beta - \lambda_\alpha)s}
					s^{k-1} \mathrm ds
					= \left( \frac{\lambda_\beta}{\lambda_\beta-\lambda_\alpha} \right)^k,
				\end{align}
				if $\lambda_\beta>\lambda_\alpha$.
				Analogously we have that
				\begin{align}
					\mathbb{V}\text{ar} Y_a(\tau_k^1) & = \frac{\lambda_\beta^k}{(k-1)!} \int_0^\infty
					e^{\lambda_\alpha s} (1+e^{\lambda_\alpha s})
					s^{k-1} e^{-\lambda_\beta s} \mathrm ds
					= \left( \frac{\lambda_\beta}{\lambda\beta-\lambda_\alpha} \right)^k
					+ \left( \frac{\lambda_\beta}{\lambda_\beta-2\lambda_\alpha} \right)^k,
					\qquad \lambda_\beta > 2 \lambda_\alpha.
				\end{align}
			\end{Remark}

		\subsection{Composition of Poisson processes with the inverse of an independent fractional
			linear birth process}
			
			Let $Y_\beta^\nu(t)$, $t>0$, be a fractional linear pure birth process with rate $\lambda_\beta>0$,
			studied in \citet{polbir}. From \citet{cah}, the distribution of
			\begin{align}
				\phi_k^\nu = \inf (t \colon Y_\beta^\nu(t) = k),
			\end{align}
			is obtained and has the following probability density:
			\begin{align}
				\label{123}
				\text{Pr} \{ \phi_k^\nu \in \mathrm dt \} / \mathrm dt 
				& = \sum_{m=1}^k \sum_{l=1}^m \binom{m-1}{l-1} (-1)^{l-1} \lambda_\beta l t^{\nu-1}
				E_{\nu,\nu}(-\lambda_\beta l t^\nu) \\
				& = \sum_{l=1}^k (-1)^{l-1} \lambda_\beta l t^{\nu-1} E_{\nu,\nu}(-\lambda_\beta l t^\nu)
				\sum_{m=l}^k \binom{m-1}{l-1} \notag \\
				& = \sum_{l=1}^k \binom{k}{l} (-1)^{l-1} \lambda_\beta l t^{\nu-1} E_{\nu,\nu}(-\lambda_\beta l t^\nu)
				\notag \\
				& = \sum_{l=1}^k \binom{k}{l} (-1)^l
				\frac{\mathrm d}{\mathrm dt} E_{\nu,1}(-\lambda_\beta l t^\nu), \qquad t>0, \: \nu \in (0,1]. \notag
			\end{align}
			The relation $\sum_{m=l}^k \binom{m-1}{l-1} = \binom{k}{l}$,
			used in the second step of \eqref{123}, can be proved as follows:
			\begin{align}
				\sum_{m=l}^k \binom{m-1}{l-1}
				& = 1 + l + \frac{l(l+1)}{2} + \frac{l(l+1)(l+2)}{2\cdot 3} 
				+ \dots + \frac{l(l+1)\dots(k-1)}{2\cdot 3 \cdots (k-l)} \\
				& = (l+1)\left[ 1+\frac{l}{2} + \frac{l(l+2)}{2\cdot 3} + \dots +
				\frac{l(l+2)\dots(k-1)}{2\cdot 3 \cdots (k-l)} \right] \notag \\
				& = \frac{(l+1)(l+2)}{2} \left[ 1+ \frac{l}{3} + \dots + \frac{l(l+3)\dots(k-1)}{
				3\cdot 4 \cdots (k-l)} \right] \notag \\
				& = \frac{(l+1)(l+2)(l+3)\dots(k-1)}{(k-l-1)!} \left[ 1+\frac{l}{k-l} \right]
				= \binom{k}{l}. \notag				
			\end{align}
			The distribution of $N_\alpha(\phi_k^\nu)$ therefore becomes
			\begin{align}
				\text{Pr} (N_\alpha(\phi_k^\nu)=r)
				& = \int_0^\infty \frac{e^{-\lambda_\alpha s}(\lambda_\alpha s)^r}{r!}
				\sum_{l=1}^k \binom{k}{l} (-1)^{l-1} \lambda_\beta l s^{\nu-1} E_{\nu,\nu} (-\lambda_\beta l s^\nu)
				\mathrm ds \\
				& = \frac{1}{r!} \sum_{l=1}^k
				\binom{k}{l} (-1)^l
				\sum_{n=0}^\infty \left( -\frac{\lambda_\beta l}{\lambda_\alpha^\nu} \right)^{n+1}
				\frac{\Gamma(\nu(n+1)+r)}{\Gamma(\nu(n+1))} \notag.
			\end{align}
			The probability generating function of $N_\alpha(\phi_k^\nu)$ can be written in a neat form as
			\begin{align}			
				\mathbb{E}u^{N_\alpha(\phi_k^\nu)} & = \sum_{r=0}^\infty u^r
				\int_0^\infty e^{-\lambda_\alpha s} \frac{(\lambda_\alpha s)^r}{r!}
				\sum_{l=1}^k \binom{k}{l} (-1)^{l-1} \lambda_\beta l s^{\nu-1} E_{\nu,\nu} (-\lambda_\beta l s^\nu)
				\mathrm ds \\
				& = \int_0^\infty e^{-\lambda_\alpha s(1-u)}
				\sum_{l=1}^k \binom{k}{l} (-1)^{l-1} \lambda_\beta l s^{\nu-1} E_{\nu,\nu} (-\lambda_\beta l s^\nu)
				\mathrm ds \notag \\
				& = \sum_{l=1}^k \binom{k}{l} (-1)^{l-1} \frac{\lambda_\beta l}{\left[ \lambda_\alpha
				(1-u) \right]^\nu + \lambda_\beta l}
				= k \sum_{l=0}^{k-1} \binom{k-1}{l} (-1)^l \frac{1}{\frac{\lambda_\alpha^\nu(1-u)^\nu}{\lambda_\beta}
				+1+l} \notag \\
				& = \frac{k!}{\left(\frac{\lambda_\alpha^\nu(1-u)^\nu}{\lambda_\beta}+1\right)
				\left(\frac{\lambda_\alpha^\nu(1-u)^\nu}{\lambda_\beta}+2\right)\cdots
				\left(\frac{\lambda_\alpha^\nu(1-u)^\nu}{\lambda_\beta}+k\right)} \notag \\
				& = \frac{k!\Gamma\left( \frac{\lambda_\alpha^\nu(1-u)^\nu}{\lambda_\beta} +1 \right)}{
				\Gamma\left( \frac{\lambda_\alpha^\nu(1-u)^\nu}{\lambda_\beta} +1+k \right)}
				= k \cdot \text{Beta}
				\left( k, \frac{\lambda_\alpha^\nu(1-u)^\nu}{\lambda_\beta} +1 \right), \qquad |u|<1.
				\notag
			\end{align}

	\section{Poisson random products and Poisson random continued fractions}	
			
		\subsection{Multiplicative compound Poisson process}
	
			In this section we consider a multiplicative compound Poisson process (denoted here $\pi$-compound Poisson
			process), defined as
			\begin{align}
				N_\pi(t) = \prod_{j=1}^{N(t)} X_j, \qquad t>0,
			\end{align}
			where the $X_j$s are i.i.d.\ random variables and $N(t)$, $t>0$, is a homogeneous Poisson process with
			rate $\lambda>0$.
			We start by calculating the Mellin transform of $N_\pi(t)$.
			\begin{align}
				\label{treccani}
				\mathbb{E} \left[ N_\pi(t) \right]^{\eta -1} = \sum_{k=0}^\infty \left[ \mathbb{E} X^{\eta -1} \right]^k
				\frac{(\lambda t)^k}{k!} e^{-\lambda t}
				= e^{\lambda t (\mathbb{E}X^{\eta -1} -1)}.
			\end{align}
			The relation \eqref{treccani} can be rewritten as
			\begin{align}
				\mathbb{E}e^{(\log N_\pi(t))(\eta -1)} = \mathbb{E} e^{(\eta -1)\sum_{j=1}^{N(t)}\log X_j}
				= e^{i\beta \sum_{j=1}^{N(t)}\log X_j} = e^{\lambda t\left( \mathbb{E}X^{i\beta} \right)}.
			\end{align}
			For the non-negative random variables $X_j$s, the random sum $\sum_j=1^{N(t)}\log X_j$ can be
			reduced to a Poisson random product for the random variables $X_j$s possessing Mellin transform at point
			$\eta=i\beta +1$,
			\begin{align}
				\mathbb{E}X^{\eta -1} = \mathbb{E} X^{i \beta}.
			\end{align}
			
			We give the explicit form of the covariance function in the next theorem.
			
			\begin{Theorem}
				For $0<s<t$, the covariance of the random product $N_\pi(t)$, $t>0$, reads
				\begin{align}
					\label{bianca}
					\mathbb{C}\text{ov} (N_\pi(t), N_\pi(s)) & =
					e^{\lambda t (\mathbb{E}X-1)} \left[ e^{\lambda s\mathbb{E}\left[ X(X-1) \right]}
					- e^{\lambda s (\mathbb{E}X-1)} \right]
					= e^{\lambda t(\mathbb{E}X-1)} \int_{s\mathbb{E}(X-1)}^{s\mathbb{E}X(X-1)} \lambda e^{\lambda w}
					\mathrm dw.
				\end{align}
				\begin{proof}
					\begin{align}
						\mathbb{E}\left[ \prod_{j=1}^{N(t)} X_j \cdot \prod_{r=1}^{N(s)} X_r \right] & =
						\mathbb{E} \left[ \prod_{j=1}^{N(s)}X_j \cdot \prod_{l=N(s)+1}^{N(t)}X_l
						\cdot \prod_{r=1}^{N(s)}X_r \right] \\
						& = \sum_{m=0}^\infty \sum_{n=m}^\infty \left[ \mathbb{E}X^2 \right]^m
						\left[ \mathbb{E}X \right]^{n-m} \text{Pr} \{ N(s)=m, N(t)=n \} \notag \\
						& = \sum_{m=0}^\infty \sum_{n=m}^\infty \left[ \mathbb{E}X^2 \right]^m
						\left[ \mathbb{E}X \right]^{n-m} \text{Pr} \{ N(s)=m \} \text{Pr} \{ N(t-s)=n-m \} \notag \\
						& = \sum_{m=0}^\infty \left[ \mathbb{E}X^2 \right]^m \frac{e^{-\lambda s}(\lambda s)^m}{m!}
						\sum_{r=0}^\infty \left[ \mathbb{E}X \right]^r \frac{e^{-\lambda(t-s)}(\lambda(t-s))^r}{r!}
						\notag \\
						& = e^{\lambda s(\mathbb{E}X^2-1)} e^{\lambda(t-s)(\mathbb{E}X-1)}. \notag
					\end{align}
					Therefore
					\begin{align}
						\mathbb{C}\text{ov}\left[ \prod_{j=1}^{N(t)} X_j, \prod_{r=1}^{N(s)} X_r \right] & =
						e^{\lambda s(\mathbb{E}X^2-1)} e^{\lambda(t-s)(\mathbb{E}X-1)}
						- e^{\lambda t(\mathbb{E}X-1)} e^{\lambda s(\mathbb{E}X-1)} \\
						& = e^{\lambda t (\mathbb{E}X-1)} \left[ e^{\lambda s\mathbb{E}\left[ X(X-1) \right]}
						- e^{\lambda s (\mathbb{E}X-1)} \right]. \notag
					\end{align}		
				\end{proof}
			\end{Theorem}
			
			\begin{Remark}
				Formula \eqref{bianca} shows that the process $N_\pi(t)$ is positively correlated.
				
				As a consequence of the previous calculations we have that
				\begin{align}
					\mathbb{E} N_\pi(t) = e^{\lambda t (\mathbb{E}X -1)}, \qquad
					\mathbb{E} \left[ N_\pi (t) \right]^2 = e^{\lambda t (\mathbb{E}X^2 -1)},
				\end{align}
				and
				\begin{align}
					\mathbb{V}\text{ar} N_\pi (t) & = e^{\lambda t (\mathbb{E}X^2 -1)} - e^{2\lambda t(\mathbb{E}X -1)}
					= e^{-\lambda t(1-\mathbb{E}X^2)} \left[ 1-e^{-\lambda t \mathbb{E}(X-1)^2}
					\right].
				\end{align}
				
				For $X \sim N(0,1)$, the covariance function of $N_\pi(t)$ takes the form
				\begin{align}
					\mathbb{C}\text{ov} N_\pi(t) & = 2 \sinh \left[\lambda \min (s,t)\right]
					e^{-\lambda \min(s,t)}
					= 2 \sinh \left[ \lambda \mathbb{C}\text{ov} (N(t),N(s)) \right]
					e^{-\lambda \mathbb{C}\text{ov}(N(t),N(s))}.
				\end{align}
			\end{Remark}

			\begin{Remark}
				If the random variables $X_j$, $j \geq 1$, are positively skewed stable with index $\nu \in (0,1)$,
				we are able to give an explicit form of the Mellin transform \eqref{treccani}.
				
				Since
				\begin{align}
					\mathbb{E}e^{-\mu X} = e^{-\mu^\nu}, \qquad \mu > 0, \: 0<\nu<1,
				\end{align}
				we have that the characteristic function of $X$ reads
				\begin{align}
					\mathbb{E}e^{i\beta X} = e^{-(i\beta)^\nu} = e^{-|\beta|^\nu e^{-\frac{i\pi \nu}{2}
					\sgn \beta}} = e^{|\beta|^\nu} 	\left[ \cos \frac{\pi \nu}{2}
					\left( 1-\sgn \beta \tan \frac{\pi \nu}{2} \right) \right].
				\end{align}
				Some manipulations as shown in \citet{mirko} prove that
				\begin{align}
					\mathbb{E}X^{\eta -1} = \frac{1}{\nu} \Gamma \left( \frac{1-\eta}{\nu} \right)
					\frac{1}{\Gamma(1-\eta)}.
				\end{align}
				This permits us to conclude that
				\begin{align}
					\mathbb{E} \left[ N_\pi(t) \right]^{\eta -1}
					& = e^{\lambda t \left( \frac{1}{\nu} \Gamma \left( \frac{1-\eta}{\nu} \right)
					\frac{1}{\Gamma(1-\eta)} -1 \right)}. 
				\end{align}
			\end{Remark}			

			\begin{Remark}
				When $X_j$ are i.i.d\ Bernoulli random variables of parameter $p$,
				we have that the fractional moments of the
				compound process can be written as
				\begin{align}
					\mathbb{E}\left[ N_\pi(t) \right]^\eta = e^{-\lambda t (1-p)},
				\end{align}
				which do not depend on $\eta$. It follows that the mean value and the variance are
				\begin{align}
					\mathbb{E} N_\pi(t) = e^{-\lambda t (1-p)}, \qquad \mathbb{V} N_\pi (t) = e^{-\lambda t(1-p)}
					\left( 1-e^{-\lambda t(1-p)}\right).
				\end{align}
				Note how the mean value and the variance formally coincide with those of a linear pure death process
				with a single progenitor.
			\end{Remark}
			
			\begin{Remark}[General case]
				Consider an infinitely divisible random variable $Y$ in the sense of Mellin (or log-infinitely
				divisible), thus
				decomposable in product of
				i.i.d.\ random variables $\zeta_j$. For $Y$ we have that
				\begin{align}
					\mathbb{E}Y^{\eta-1} = \left[ \mathbb{E}\zeta^{\eta-1} \right]^k.
				\end{align}
				The Mellin transform of the random product $\prod_{j=1}^{N(t)}\zeta_j$ is therefore
				\begin{align}
					\label{belalm}
					\mathbb{E}\left[ \prod_{j=1}^{N(t)}\zeta_j \right]^{\eta-1} & =
					\sum_{k} \left[ \mathbb{E} \zeta^{\eta-1} \right]^k \frac{e^{\lambda_\beta t}
					(\lambda_\beta t)^k}{k!}
					= e^{-\lambda_\beta t} e^{\lambda_\beta t \mathbb{E}\zeta^{\eta-1}}
					= e^{-\lambda_\beta t\left[ 1-\mathbb{E} \zeta^{\eta-1} \right]}.
				\end{align}
				Let now	$\Theta \colon \mathbb{N} \rightarrow \mathbb{N}$ such that for each
				$r \in \mathbb{N}$, $\Theta(r) = \prod_{j=1}^r \xi_j$. If the random
				variables $\xi_j$s take integer values, then
				\begin{align}
					\label{water}
					\mathbb{E}\left[ \Theta(N(t)) \right]^{\eta-1} & = \sum_{m=0}^\infty
					\sum_{r=0}^\infty m^{\eta-1} \text{Pr} \{ \Theta(r) = m \} \text{Pr} \{ N(t) = r \} \\
					& = \sum_{r=0}^\infty \left( \mathbb{E} \xi^{\eta-1} \right)^r \text{Pr} \{ N(t) = r \}
					= e^{-\lambda_\beta t + \lambda_\beta \mathbb{E}\xi^{\eta -1}}. \notag
				\end{align}
				If $\Theta(r)$ is absolutely continuous the calculations follow in the same way and arrive
				at the Mellin transform \eqref{water}:
				\begin{align}
					\label{fire}
					\mathbb{E}\left[ \Theta(N(t)) \right]^{\eta-1} & =
					\sum_{r=0}^\infty \int_0^\infty x^{\eta-1} \text{Pr} \{ \Theta(r) \in \mathrm dx \}
					\text{Pr} \{ N(t) = r \} \\
					& = \sum_{r=0}^\infty \left( \mathbb{E} \xi^{\eta-1} \right)^r \text{Pr} \{ N(t) = r \}
					= e^{-\lambda_\beta t + \lambda_\beta \mathbb{E}\xi^{\eta -1}}. \notag
				\end{align}
				In conclusion we have that
				\begin{align}
					\Theta(N(t)) \overset{\text{d}}{=} \prod_{j=1}^{N(t)} \xi_j .
				\end{align}
			\end{Remark}
				
		\subsection{Continued fractions of Cauchy random variables with Poisson distributed levels}
		
			We consider in this section the random variables defined as
			\begin{align}
				[ X_1;X_2,\dots,X_{N(t)} ] = X_1 + \frac{1}{X_2 + \frac{1}{\ddots + X_{N(t)-1} + \frac{1}{X_{N(t)}}}},
			\end{align}
			where $X_j$, $j\geq 1$, are independent Cauchy random variables with scale parameter equal to unity
			and location parameter equal to zero. We will write $X \sim C(0,1)$.
			Furthermore $N(t)$, $t>0$, is a homogeneous Poisson process independent of the Cauchy random variables
			$X_j$.
			For the convenience of the reader we note that
			\begin{align}
				[X_1] & = X_1, \\
				[X_1;X_2] & = X_1 + \frac{1}{X_2},  \\
				[X_1;X_2,X_3] & = X_1 +\frac{1}{X_2 + \frac{1}{X_3}}.
			\end{align}

			The standard Cauchy random variable has the remarkable property that $X \sim 1/X$, and this is the reason for
			which continued fractions can be treated when Cauchy random variables are involved (\citet{camma}).
			
			For our analysis, we need the following result.
			\begin{Lemma}
				\label{lem}
				For a Cauchy random variable $C(a,b)$, $a \in \mathbb{R}$, $b \in \mathbb{R}^+$, the
				random variable $1/C(a,b) \sim C(a/(a^2+b^2), b/(a^2+b^2))$. In our case, $a=0$ and therefore
				$1/C(0,b) \sim C(a,1/b)$.
			\end{Lemma}
			
			Our first result is stated in the next theorem.
			
			\begin{Theorem}
				\label{teo98}
				The $n$th level fraction
				\begin{align}
					[ X_1;X_2,\dots,X_n ] = X_1 + \frac{1}{X_2 + \frac{1}{\ddots + X_{n-1} + \frac{1}{X_n } } } ,
				\end{align}
				has Cauchy distribution with scale parameter $b_n=F_{n+1}/F_n$, where $F_n$ are the Fibonacci numbers.
				\begin{proof}
					We proceed by induction.
					\begin{align}
						[X_1;X_2] = X_1 + \frac{1}{X_2} \sim C(0,2).
					\end{align}
					In view of Lemma \ref{lem}, we have that
					\begin{align}
						[X_1;X_2,X_3] = X_1 + \frac{1}{X_2 + \frac{1}{X_3}} \sim C(0,3/2).
					\end{align}
					Furthermore,
					\begin{align}
						[X_1;X_2,X_3,X_4] = X_1 + \frac{1}{[X_2;X_3,X_4]} \sim C(0,5/3).
					\end{align}
					In general we have that
					\begin{align}
						[X_1;X_2,\dots,X_n] & = X_1 + \frac{1}{[X_2,X_3,\dots,X_n]} 
						= X_1 + \frac{1}{C(0,F_{n-1}/F_n)} \\
						& = X_1 + C(0,F_n/F_{n-1})
						= C(0,1+F_n/F_{n-1})
						= C(0, F_{n+1}/F_n), \notag
					\end{align}
					and in the last step we took into account the definition of Fibonacci numbers.
				\end{proof}
			\end{Theorem}
			
			\begin{Remark}
				The Fibonacci numbers can be written in terms of the golden ration $\phi = (1+\sqrt{5})/2$ as
				\begin{align}
					F_n = \frac{\phi^n - (1-\phi)^n}{\sqrt{5}}.
				\end{align}
				Therefore
				\begin{align}
					\frac{F_{n+1}}{F_n} & = \frac{\phi^{n+1}-(1-\phi)^{n+1}}{\phi^n + (1-\phi)^n} 
					= \phi \frac{1-\left( \frac{1-\phi}{\phi} \right)^{n+1}}{1-\left( \frac{1-\phi}{\phi} \right)^n}
					\rightarrow_{n \rightarrow \infty} \phi.
				\end{align}
				This means that $[X_1;X_2,\dots,X_n] \overset{\text{d}}{\rightarrow} C(0,\phi)$.
			\end{Remark}
			
			\begin{Remark}
				From the analysis above, we infer that $[X_1;X_2,\dots,X_{N(t)}]$, $t>0$, is a process and,
				for each $t$, possesses distribution equal to
				\begin{align}
					\text{Pr} \{ [X_1;X_2,\dots,X_{N(t)}])
					\in \mathrm dx \} /\mathrm dx= \sum_n \frac{1}{\pi} \frac{F_{n+1}/F_n}{x^2+
					(F_{n+1}/F_n)^2} e^{-\lambda t} \frac{(\lambda t)^n}{n!}.
				\end{align}		
			\end{Remark}
			We give now an alternative representation to the process $[X_1;X_2,\dots,X_{N(t)}]$, $t>0$.
			\begin{Theorem}
				The characteristic function of the random continued fraction $[X_1;X_2,\dots,X_{N(t)}]$
				reads
				\begin{align}
					\label{baffioni}
					\mathbb{E}e^{i\beta [X_1;X_2,\dots,X_{N(t)}]} =
					e^{-|\beta|\phi} \sum_{n=0}^\infty \prod_{j=1}^\infty e^{-|\beta|
					\sqrt{5} \left( \frac{1-\phi}{\phi} \right)^{nj}} \text{Pr} \{ N(t)=n \}.
				\end{align}
				\begin{proof}
					In view of Theorem \ref{teo98}, we have that
					\begin{align}
						\label{asterisque}
						\mathbb{E} e^{i\beta [X_1;X_2,\dots,X_{N(t)}]} & =
						\sum_{n=0}^\infty \mathbb{E} e^{i\beta [X_1;X_2,\dots,X_n]} \text{Pr} \{ N(t)=n \} \\
						& = \sum_{n=0}^\infty e^{-|\beta| \frac{F_{n+1}}{F_n}} \text{Pr} \{ N(t) = n \}. \notag
					\end{align}
					Since
					\begin{align}
						\frac{F_{n+1}}{F_n} & = \phi \frac{1-\left( \frac{1-\phi}{\phi} \right)^{n+1}}{1-\left(
						\frac{1-\phi}{\phi} \right)^n} \\
						& = \phi \left\{ 1-\left( \frac{1-\phi}{\phi} \right)^{n+1} \right\}
						\sum_{j=0}^\infty \left( \frac{1-\phi}{\phi} \right)^{nj} \notag \\
						& = \phi \left\{ \sum_{j=0}^\infty \left( \frac{1-\phi}{\phi} \right)^{nj}
						-\frac{1-\phi}{\phi} \sum_{j=0}^\infty \left( \frac{1-\phi}{\phi} \right)^{n(j+1)} \right\}
						\notag \\
						& = \phi \left\{ 1 + \left( 1-\frac{1-\phi}{\phi} \right)
						\sum_{j=1}^\infty \left( \frac{1-\phi}{\phi} \right)^{nj} \right\} \notag \\
						& = \phi + \sqrt{5} \sum_{j=0}^\infty \left( \frac{1-\phi}{\phi} \right)^{nj}, \notag
					\end{align}
					we have that
					\begin{align}
						\mathbb{E} e^{i\beta [X_1;X_2,\dots,X_n]} & =
						\sum_{n=0}^\infty e^{-|\beta| \left[ \phi + \sqrt{5}
						\sum_{j=1}^\infty \left( \frac{1-\phi}{\phi} \right)^{nj} \right]}
						\text{Pr} \{ N(t)=n \} \\
						& = e^{-|\beta|\phi} \sum_{n=0}^\infty \prod_{j=1}^\infty e^{-|\beta|\sqrt{5}
						\left( \frac{1-\phi}{\phi} \right)^{nj}} \text{Pr} \{ N(t)=n \}. \notag
					\end{align}
				\end{proof}
			\end{Theorem}
			
			\begin{Remark}
				From \eqref{baffioni} we can extract the following equality in distribution:
				\begin{align}
					\label{ago}
					[X_1;X_2,\dots,X_{N(t)}] \overset{\text{d}}{=} C(0,\phi) + \sqrt{5}
					\sum_{j=1}^\infty C_j\left( 0,\left( \frac{1-\phi}{\phi} \right)^{N(t)} \right).
				\end{align}
				The second term in \eqref{ago} represents the effect of randomisation of the continued fraction.
			\end{Remark}
			
			\begin{Remark}
				The above analysis suggests an alternative representation of the random continued
				fraction as
				\begin{align}
					\label{stilo}
					[X_1;X_2,\dots,X_{N(t)}] \overset{\text{d}}{=} \sum_{j=1}^{F_{N(t)+1}}
					Y_{j,N(t)}, 
				\end{align}
				where the $Y_{j,N(t)}$ are independent Cauchy random variables with scale parameter equal to
				$1/F_{N(t)}$. Clearly, $F_n$ are the Fibonacci numbers.
				
				The equality \eqref{stilo} can be ascertained by writing the characteristic function
				as follows:
				\begin{align}
					\mathbb{E} e^{i\beta [X_1;X_2, \dots, X_{N(t)}]} & =
					\mathbb{E} e^{i\beta \sum_{j=1}^{F_{N(t)+1}}Y_{j,N(t)}}
					= \mathbb{E}\left[ \left. \mathbb{E} e^{i\beta \sum_{j=1}^{F_{N(t)+1}}Y_{j,N(t)}}
					\right| N(t) \right] \\
					& = \sum_{n=0}^\infty \prod_{j=1}^{F_{n+1}} e^{-|\beta|\frac{1}{F_n}} \text{Pr} \{ N(t) = n \}
					= \sum_{n=0}^\infty e^{-|\beta|\frac{F_{n+1}}{F_n}} \text{Pr} \{ N(t)=n \}, \notag
				\end{align}
				which coincides with \eqref{asterisque}.
			\end{Remark}
			
	\bibliographystyle{abbrvnat}
	\bibliography{comp-poisson3}
	\nocite{*}

\end{document}